\documentclass{amsart}
\usepackage{eepic}
    \usepackage{a4wide} 

\newif\ifview
\viewfalse  
\ifview
    \scrollmode  
    \def\herea#1.{\marginpar{#1 \vrule width 2pt depth 1pt height 9pt \ }}
\else
   \def\herea#1.{}

\fi

\usepackage{amsthm}
\theoremstyle{remark} 
\newtheorem{thm}{xx}[section]

\def\newt#1 {\newtheorem{#1}[thm]{#1}}
\newt Definition
\newt Theorem
\newt Corollary
\newt Lemma
\newt Remark
\newt Question
\newt Notation
\newt Fact
\newt Note
\newt Convention
\newt Proposition
\newt Conclusion
\newt Naming
\newt Construction
\newt Motivation
\newtheorem{CL}[thm]{Crucial Lemma}
\newtheorem{mn}[thm]{More Notation}
\usepackage{mathrsfs}
\usepackage{amssymb}

\newcommand{\N}{{\mathbb N}}
\newcommand{\Z}{{\mathbb Z}}

\newcommand{\on}{{\upharpoonright}}
\def\xx#1^#2 {#1^{(#2)}}  
 \def\yy#1^#2{#1^{[#2]}}  
 \newcommand{\leaf}{{\bf ext}}   
 \newcommand{\ext}{\leaf} 
  
 \renewcommand{\succ}{{\bf succ}} 
  \newcommand{\T}{{ T}} 
  \renewcommand{\S}{{ S}} 

 \renewcommand{\int}{{\bf int}} 
 \renewcommand{\O}{{\mathscr O}} 
 \newcommand{\C}{{\mathscr C}}

 \newcommand{\cid}{{\C_{\rm id}}} 
 \newcommand{\id}{{\rm id}} 
 \newcommand{\wurzel}{{\bf root}} 
 \newcommand{\J}{{\mathscr J}} 
 \newcommand{\ramsey}{{\mathscr R}} 
 \newcommand{\F}{{\mathscr{F}}} 
\newcommand{\G}{{\mathscr{G}}} 
 
 \newcommand{\kk}{\hbox to 0pt{\hss $\scriptscriptstyle k\in
 \ext(r)$\hss}} 
 \newcommand{\kkn}{\hbox to 0pt{\hss $\scriptscriptstyle k\in
 \ext(R_n)$\hss}}  
\newcommand{\Pol}{{\rm Pol}}
\newcommand{\I}{{\mathscr I}}
\newcommand{\prs}{\preccurlyeq^*_{\vec s}}
\newcommand{\aps}{\approx^*_{\vec s}}
  
 \def\oo#1 {{\O^{(#1)}}} 
 
\begin{document} 
  
 \title{Clones from Creatures} 
 \author{Martin Goldstern} 
 \address{Vienna University of Technology} 
 \thanks{The first author is grateful to the Math department of 
  Rutgers University, New Jersey, for  
 their hospitality during a visit in September 2002} 
 \author{Saharon Shelah} 
 \address{Hebrew University of Jerusalem and Rutgers University, NJ} 
  
 \thanks{The second author's research  was supported by 
   the United States-Israel Binational Science Foundation. 
 Publication 808.}

 \begin{abstract}  We show that (consistently) there is a clone $\C$  
 on a countable set such that the interval of clones above $\C$ is 
 linearly ordered and has no coatoms. 
 \end{abstract} 
  
\date{Feb 24, 2003}
  
 \maketitle

\ifview  \markboth{XXX}{XXX}\fi
  
 \setcounter{section}{-1} 
 \section{Introduction }  
  
 A clone on a set $X$ is a set of finitary operations $f:X^n\to X$ 
 which contains all the projections and is closed under 
 composition. 
 (Alternatively, $\C$ is a clone {if} $\C$ is the set of term functions of 
 some universal  algebra over~$X$.)  
  
  The family of all clones forms a complete algebraic 
 lattice $Cl(X)$ with greatest element $\O = \bigcup_{n=1}^\infty \oo n $,  
 where $\oo n = X^{X^n} $ is the set of all $n$-ary operations on~$X$.
   (In this paper, the underlying set $X$  will always be the set 
 $\N = \{0,1,2,\ldots \} $ of natural numbers.) 
  
   The coatoms of this lattice $Cl(X) $ are called 
 ``precomplete clones'' or ``maximal clones'' on~$X$.  
  
 For singleton sets $X$ the lattice $Cl(X)$ is trivial; for $|X|=2$
 the lattice $Cl(X)$ is countable, and well understood (``Post's
 lattice'').   
 For $3\le |X| < \aleph_0$, $Cl(X)$ has uncountably many elements.  
 Many results for clones on finite sets 
  can be found in  
 \cite{Szendrei:1986}.    In particular, there is an explicit
 description of all  
 (finitely many) precomplete clones on a given finite set 
  (\cite{Rosenberg:1970}, see also   
   \cite{Quackenbush:1971} and \cite{Buevich:1996});
 this description also includes a decision procedure
 for the membership problem for each of these clones. 
It is also known 
 that every clone $\C\not=\O$ is contained in a precomplete 
 clone, that is: the clone lattice  $Cl(X)$  on any finite set  
 $X$ is {\em dually atomic}. 
  (This gives an explicit criterion for deciding whether a given 
 set of functions generates all of~$\O$: just check {if} it is contained 
 in one of the precomplete clones.)  
  
 Fewer results are known about the lattice of clones on an infinite 
 set: \cite{Gavrilov:1965} investigated the interval of clones above 
 the clone of unary functions on a countable set,  \cite{GoSh:737} did 
 this also for uncountable sets.  \cite{Heindorf:2002} classified the 
 countably many precomplete clones on a countable set that contain all 
 bijections.  
   \cite{Rosenberg:1976} showed that 
 there are always $2^{2^\kappa}$ precomplete clones on a set of 
 infinite cardinality~$\kappa$, and \cite{Rosenberg:1974} gave specific 
 examples of such precomplete clones. \cite {Machida+Rosenberg:1992} 
 investigated minimal clones.  
  
  \cite {Rosenberg+Schweigert:1982} 
 investigated  
 ``local'' clones on infinite sets (clones that are closed sets  in the 
 product topology).  It  is easy to see that the lattice of local 
 clones is far from being dually atomic (\cite{GoSh:747}).

Already Gavrilov in \cite[page 22/23]{Gavrilov:1959} asked whether the lattice
 of all clones on a countable set is also dually atomic, since a
 positive answer would be an important component for a completeness
 criterion, as remarked above.  This question is also listed as
 problem P8 in \cite[page 91]{PK:1979}, and has been open until now.

 We will show here (assuming the continuum hypothesis, CH)
 the following:  
 \begin{Theorem}[CH]   
 \herea2. 
  \label{thm} 
 The lattice of clones on a countably 
 infinite set is {\bf not dually atomic}, i.e., there is a clone 
 $\C\not=\O$ which is not contained in any precomplete clone.  
 \end{Theorem} 
We also remark that the full strength of CH ist not needed for this theorem. 

 The clone $\C_U$ that we construct has the additional feature that we can 
 give a good description of the interval~$[\C_U,\O]$.  (See \ref{cantor1}(2) 
and \ref{done}(d).)  In particular, 
 it will be a {\em linear} order without a penultimate element,
in which every countable set  has an
 upper bound. 
 
 All clones that we consider will be in the interval~$[\cid,\O]$, where 
 $\cid$ is the clone  
 of all functions which are  bounded by the $\max$ function. The 
 clones in the interval $[\cid,\O]$ have the property that they are 
 determined by their unary part.   Moreover, the map that assigns to 
 each clone its unary part is a lattice  
  isomorphism between $[\cid, \O]$ and the set of all those monoids 
 $\subseteq \N^\N$ 
 which are lattice ideals (in the product order, see
 definition~\ref{ideal}).

 Thus, our theorem can be reformulated as follows:  
 \begin{Theorem}[CH] 
 \herea4. 
  The lattice  of all submonoids of~$\N ^\N  $ which are ideals   
 is not dually atomic.  
 \end{Theorem}

The method behind our proof is ``forcing with normed creatures'', a
set-theoretic  
construction originating in the second author's paper \cite{Sh:207}.   
 The book \cite{RoSh:470}, an encyclopedia of such creatures, may be 
 useful for constructing variants of our clone $\C_U$ to get clone 
 intervals with prescribed properties;  however, for the purposes of 
 this paper the connection with forcing machinery is sufficiently 
 shallow to allow us to be self-contained.  
  
 In particular, no knowledge  
 of set theory 
  is required for our theorem, except for the last section, where 
a basic understanding of CH and transfinite induction up to  $\omega_1$ is needed. 
 Most of our constructions deal with
 finite structures, or with countable sequences of finite structures.

 \section{Clones defined by growth conditions}
  
\label{section.reduce}

 \begin{Notation}\label{notation5} 
 \herea5. 
 \begin{enumerate} 
 \item $\oo 1  =  \N ^\N$ is the monoid of all functions from $\N$ to
  $\N$ (the operation is composition of functions.)   
  For $k\ge 1$,   
 $\oo k $ is the set of all functions from $\,  \N ^k$ to~$\N$, and   
  $\O = \bigcup_{k=1}^\infty  \oo k $.  
 \item For $f,g\in \N ^\N$, we write  $f\le g$ for $\forall n: f(n)\le g(n)$.  
 	 
 \item We write $f\le^* g$ iff $f(n) \le g(n) $ holds almost all (i.e., 
for all but finitely many)~$n$.   In general we use the superscript
 $\ {}^*\ $ 
or the keyword ``almost'' to
 indicate that finitely many   
 exceptions are allowed.  
 \item $\id:\N\to \N$ is the identity function. 
 \item For $f\in \N   ^\N$, we write  $\xx f^n $  for the $n$-fold
 composition of~$f$ with itself ($\xx f^0 = \id$, $\xx f^1 = f$.)  
 \item $\max_k$ is the $k$-ary maximum function.  Usually we  just 
write~$\max$. 
 \item  \label{growth}
   A {\em growth function} is a (not necessarily strictly) 
  increasing function from $\N$ to~$\N$  satisfying $f(n)>n$ for all~$n$.  
We write $\G$ for the set of all growth functions.

 \item \label{gbar}
For any function $g:\N^k\to \N$ we let $\bar g:\N\to \N$  be
 defined by   
 $$ \bar g (n) := \max \{ g(x_1,\ldots, x_k): x_1,\ldots, x_k\le n \}, 
 $$ and we let $\hat g(n) = \max\{n, \bar g(n)\}$.  So $\hat g +1$ is a 
 growth function. 
 \item If $f$ is $k$-ary, $g_1, \ldots, g_k$ are $n$-ary,  
    we write $ f(g_1,\ldots, g_k) $ for the function that maps $\vec x =  
       (x_1,\ldots, x_n) $ to~$f(g_1(\vec x), \ldots, g_k(\vec x)) $. 
      For example,  $\max(g,h)$ is the pointwise maximum  of~$f$  and~$g$. 
 \item $\cid:= \{ g\in \O: g \le \max \} = 
\{ g\in \O: \bar g \le \id \} $. 
 \end{enumerate} 
 \end{Notation}

 It is clear that $\cid$ is a clone.    We will only consider clones that 
 include~$\cid$.  
  
 Note that the function  $\bar g$ is always  increasing (not necessarily 
 strictly).  $\bar g$ measures the ``growth'' of~$g$. 
  
 \begin{mn} 
 The following symbols are
 collected here   only for easier reference:  
 \begin{itemize} 
 \item Relations between functions: $f\le g$, $f\le^*g$, $f= ^* g$:
 see~\ref{notation5}.    
\item Growth functions,~$\G$: \ref{notation5}.
 \item The function $h_A$  and the relation  
   $f\le_A g$ for infinite sets~$A$:~\ref{star}.   
 \item Relations between fronts or *fronts in zoos:  
 $F \prec_s G$, $F\prec^*_s G$, $F\preccurlyeq_s G$, $F\preccurlyeq^*_s 
 G$, $F \approx^*_s G$: see~\ref{prec}. 
 \item More relations between fronts:  
 $F + n \preccurlyeq_s^*  G$, $F + n \approx _s^*  G$, 
 $F + \infty  \preccurlyeq_s^*  G$. Again see~\ref{prec}.  
 \item Relations between  growth  
 functions $f$ and $g$ that are gauged by a zoo~$s$:  
  $f \prec_s g$, $f \prec_s^* g$, $f\preccurlyeq_s^* g$, $f
 \approx_s^* g$.   
 See~\ref{precf}. 
 \item The *fronts~$F(s,f)$:~\ref{front} 
 \item Relations between zoos:  
 $t= ^* s $:~\ref{defstar},
 $t\le s$,  $t \le ^* s $:~\ref{defle},  
 $t\leqq s$, $t \leqq^* s $:~\ref{defqq}.  
 \item Tree order~$\vartriangleleft$, $\trianglelefteq$:~\ref{def20}.   
 Lexicographic order $<$ on nodes in
 creatures:~\ref{def28}, and on zoos:~\ref{def43}. 
 Direct lexicographic successor $\lessdot$:~\ref{def43}.  
 \end{itemize} 
 \end{mn}

  \subsection*{From clones to ideal monoids} 

  We first show that above $\cid$ we can restrict our attention to unary functions. 
  
 \begin{Lemma}\label{lem7} 
 \herea7. 
 Let $\C$ be a clone with $\cid \subseteq \C$.  Then: 
 \begin{enumerate} 
 \item $\C$ is downward closed: If $f\le g$, and $g\in \C$, then~$f\in
 \C$.   
 \item For all~$g:\N^k\to \N$:  $g\in \C$ iff  $\bar g\in \C$  
 iff $\hat g\in \C$.  (See \ref{notation5}(\ref{gbar}).)
 \item  
  If the successor function $\id + 1 $ (mapping  each $x\in \N$ to~$ x+1 $)
 is in~$\C$, then we also have $g\in
  \C $ iff $\hat g + 1 \in \C$.    
 \end{enumerate} 
 \end{Lemma} 
 \begin{proof}

 (1):  Let $g\in \C$, $g:\N^k\to \N$, and~$f\le g$.  
 Define  a $k+1$-ary function $F$ by  
 $$ F(\vec x, n) = \min ( f(\vec x), n)$$ 
 Clearly $F\in \cid$, and $f(\vec x)  =  F(\vec x, g(\vec x))$ for all  
 $\vec x = (x_1,\ldots, x_k)$,  so~$f\in \C$.  
  
 \medskip 
  
(2): Note that $g \le \bar g(\max_k)$, and $\bar g \le \hat g$,  
  so $\hat g\in \C \Rightarrow \bar g \in\C \Rightarrow g\in \C$,
as $\C$ is downward closed.  
 The implication $\bar g\in \C\Rightarrow\hat g\in \C$ follows 
from~$\max_2\in    \cid$.  
  
 It remains to check $g\in \C\Rightarrow \bar g\in \C$: \\ 
 Assume~$g\in \C$.  
 For each $n\in \N$ choose 
 $(h_1(n),\ldots, h_k(n))\le (n,n,\ldots, n)$ such that  
 $$ g(h_1(n),\ldots, h_k(n)) =  
\max\{g(x_1,\ldots, x_k): x_1,\ldots, x_k\le n\}$$ 
  
 Then~$h_1$, \dots , $h_k\in \cid$, so $\bar g = g(h_1,\ldots, h_k)\in \C$.  
  
 (3):   The implication $\hat g \in \C \Rightarrow (\hat g + 1) \in \C$
 follows because $\id + 1 \in \C$, the converse is true because $\C $ 
is downward closed.
 \end{proof}

 \begin{Definition} 
 \herea10. 
 \label{ideal}  
 A set $M \subseteq \oo 1 $ is a {\em ideal monoid} iff $M$ is both a monoid  
 and a (lattice) ideal, i.e.:   
 \begin{enumerate} 
 \item $(M,{\circ},\id)$ is a monoid 
 \item $M$ is downward closed: $ \forall g\in M\, \forall f\le g\,\, f\in M$ 
 \item $M$ is closed under~$\max$: $\forall f_1,f_2\in M: 
 \max(f_1,f_2)\in M$.  
  \\{} 
 Using (1) and (2), this is equivalent to: $\forall f\in M: 
 \max(f,{\rm id})\in M$.  
 \end{enumerate} 
 Let $\J$ be the set of ideal monoids.  
 \end{Definition} 
  
 \begin{Proposition} 
 \herea12. 
 \begin{enumerate} 
 \item $( \J, {\subseteq })$ is a complete algebraic lattice, isomorphic 
 to the interval $[\cid,\O]$ in the clone lattice. 
 \item The map $\C \mapsto \C\cap \oo 1 $ is an isomorphism from   
 $[\cid,\O]$ onto~$\J$, with inverse $M \mapsto \{g\in \O: \bar g\in M\}$.  

 \end{enumerate} 
 \end{Proposition} 
 \begin{proof} If $\C\supseteq \C_\id$ is a clone, then $\C$ is
 downward closed by~\ref{lem7}(1) and contains the $\max$ function, so
 $\C\cap \oo 1  $ is an ideal monoid. 

Conversely, if $M$ is an ideal monoid then $\C(M):=\{g:\bar g\in M\}$
certainly contains~$\C_\id$.   To check that $\C(M)$  is closed under
composition of functions it is enough to verify 
\begin{quote}
 If $h = g(f_1, \ldots, f_k)$, and $\bar f_1,\ldots, \bar f_k, \bar
 g\in M$, then also $\bar h \in M$
\end{quote}
which follows from $\bar h \le \bar g (\max(\bar f_1,\ldots, \bar
f_k))$. 

The fact that the two maps are inverse to each other follows 
 easily from lemma~\ref{lem7}(2).   
 \end{proof} 

\begin{Remark} The isomorphism from~$\J$ onto~$[\C_\id,\O]$ can also 
be described by the map~$M\mapsto \Pol(M)$.  
\end{Remark}

So from now on we will only investigate ideal monoids instead of clones 
above~$\C_\id$.  Our aim is to
find an ideal monoid $M$ such that the interval $[M, \oo 1 ] $ of
ideal monoids (is linearly ordered and) has no coatom.

\subsection*{From ideal monoids to growth semigroups} 
The next step is mainly cosmetical; functions $f$ satisfying $f(n)=n$ 
for some $n$ are unpleasant to work with, so we want to ignore them. 

\def\itm#1 {\item[(#1)]}

Recall that $\G$ is the set of growth functions
(see \ref{notation5}(\ref{growth})).   
Note the following easy facts: 

\begin{Fact} 
\begin{enumerate}
\itm a If $f$ and $g$ are growth functions, then $f\circ g$ is a
growth function.  
\itm b For any $n\in \N$, the set $S_n:= \{f\in \oo 1
\mid \forall k\le n: f(k)\le n\}$   (sometimes written
 $\Pol(\{1,\ldots, n\})$) is a coatom in the lattice of ideal monoids. 
\itm c
The map $S \mapsto S\cap \G$, mapping an ideal monoid to its set of
growth functions, is not 1-1.
\end{enumerate}
\end{Fact}
\begin{proof}  (a) is trivial, and (b) is easy and well-known. 

For (c), note that the sets $S_n$ from (b) are all different, but all
satisfy $S_n \cap \G = \emptyset$.
\end{proof}

These observations motivate us to restrict our attention to the set of
those ideal monoids that contain the function~$x\mapsto x+1$ (the
smallest growth function).

\begin{Definition}  
A set $S \subseteq \G$ is called a growth semigroup iff $S$ is
nonempty,  closed
under composition of functions, and (in $\G$) downward closed. 
(Since $\max(f,g)\le f\circ g$ holds for all growth functions $f$ and~$g$, 
$S$ must also be closed under the $\max$ function.)

We write $\C_{\id+1}$ for the clone generated by~$\C_\id$ 
together with the successor function
$x\mapsto x+1$; thus $\C_{\id+1}\cap \oo 1 $ is the smallest ideal
monoid containing the successor function, and 
 $\C_{\id+1}\cap \G$  is the smallest  growth semigroup. 
\end{Definition}

\begin{Fact}
The map $S \mapsto S\cap \G$ is an isomorphism from the 
interval $[\C_{\id+1}\cap \oo 1 , \oo 1 ]  $ of ideal monoids onto 
the set of growth semigroups. 
\end{Fact}

\begin{proof}
   By lemma~\ref{lem7}(3). 
\end{proof}

 So the interval $[\C_{\id+1}, \O]$ in the clone lattice is isomorphic
 to the lattice of growth semigroups. From 
now on we will only investigate growth semigroups. 
Our aim is to
find a growth semigroup  $S$ such that the interval $[S, \G ] $ of
growth semigroups (is linearly ordered and) has no coatom.

\subsection*{From growth semigroups to single growth functions}

The next reduction is the most important one.  Instead of
investigating a lattice of growth  semigroups
 (or clones, or ideal monoids),  we can
reduce our analysis to the investigation of the natural 
partial quasiorder  ``$g$ generates $f$''
 of growth functions which (after factorization) turns out 
to be an upper semilattice.  The interval of clones/monoids/semigroups 
that we are interested in will be naturally isomorphic to the 
set of ideals of 
this semilattice.

\begin{Fact}\label{leqS}
Let $S$ be a growth semigroup,   $f$ and $g$ growth functions. Then the 
following are equivalent: 
\begin{enumerate}
\item  $f  $ is in the smallest growth semigroup containing $S \cup \{g\}$.
\item There is a growth  function $h \in S$ and a
natural number $k$ such that
 $f \le \xx (h \circ g)^k $ 
\item  There is a growth function $h \in S$ and a
natural number $k$ such that
 $f \le \xx (\max(h,  g))^k $ 
\item $f$ is in the clone generated by~$S \cup \C_\id \cup \{g\}$. 
\end{enumerate}
\end{Fact}
\begin{proof} 
 The equivalence of (2) and (3) follows from $\max(h,g)\le h\circ g
 \le \xx (\max(h,g))^2 $.  The implication ``(2)$\Rightarrow$(1)'' 
follows  from the closure properties of growth semigroups, and 
 ``(1)$\Rightarrow$(2)'' follows from the fact that, for 
any~$g$,  the set of growth 
functions $f$ satisfying (2) is a growth semigroup. 

The equivalence of (4) to the other conditions is left to the reader. 
\end{proof}
\begin{Definition}
If any/all of the above conditions (1), (2), (3), (4)
 are satisfied, we will write $f
\le_S g$, and we will write $\sim_S$ for the associated
equivalence relation: $f\sim_S g $ iff $f \le _S g$ and 
 $g \le _S f$.   We write $\G/{\sim_S}$ or just $\G/S$ 
for the set of equivalence classes.


We have:  
\begin{enumerate}
\itm a $\le_S$ is a partial quasiorder of functions. 
\itm b The set $S$ is the smallest $\sim_S$-equivalence class
\itm c The set $\G/{S}$ (ordered naturally) is a join-semilattice. 
\end{enumerate}
\end{Definition}

\begin{proof} 
It is clear that $\le_S$ is transitive and reflexive. Clause
 (1) in the definition of~$\le_S$
 easily implies that $S$ is the smallest
 equivalence class. 

We will write $f/S$ for the $\sim_S$-class of~$f$. 
To show clause (c) we prove that the class  $\max(f_1,f_2)/{S} $ 
is the least upper bound of the classes $f_1/{S}$ and~$f_2/{S}$.  Clearly 
$\max(f_1,f_2)/S$ is an upper bound; {if} $g/S$ is also an upper bound, then 
we have $f_1\le\xx( \max(h_1,g))^{k_1} $  and   
$f_2\le\xx( \max(h_2,g))^{k_2} $ for some $h_1,h_2\in S$, $k_1,k_2\in \N$. 

Letting  $h:=\max(h_1,h_2)\in S$ and $k:= \max(k_1,k_2)$ we have 
$\max(f_1,f_2) \le \xx( \max(h,g))^{k} $, so
$\max(f_1,f_2)\le_S g$.   
\end{proof}

\begin{Fact}  Fix a growth semigroup~$S_0$. 

Then for every growth semigroup $S \supseteq S_0$, the 
set $\{ f/S_0:  f\in S \}$  is an ideal in
the semilattice~$\G/S_0$. 

Conversely, {if} $I \subseteq \G/S_0$ is a (nonempty) semilattice ideal,
then 
$\{ f\in \G: f/S_0\in I\}$ is a growth semigroup containing~$S_0$. 

Moreover, the maps defined in the previous two paragraphs are inverses
of each other. Thus, the interval $[S_0,\G]$ in the set of growth
semigroups is naturally isomorphic to the set of ideals on~$\G/S_0$.

In particular we get:
  If $\G/S_0$ is linearly ordered, then
\begin{enumerate}
\item  the interval 
$[S_0, \G]$ in the lattice of growth semigroups corresponds 
exactly to the nonempty downward closed subsets  of~$\G/S_0$ (the
``Dedekind cuts''). 
\item $[S_0, \G]$ has a coatom iff $\G/S_0$ has a greatest element. 
\end{enumerate}

\end{Fact}

\begin{proof} 
Again this boils down to 
$\max(h,g)\le h\circ g  \le \xx (\max(h,g))^2 $.    The fact that 
semigroups are always unions of~$\sim_{S_0}$-equivalence classes 
follows from the relation~$f \circ f \sim_{S_0} f$. 


\end{proof}

\subsection*{Conclusion and goals}

So far we have shown that to find a clone which is not below a coatom
it is enough to find a growth semigroup ${S_0}$ such that the set of  growth
functions, partially quasiordered by~$\le_{S_0}$, has no maximal element.

We will construct a growth semigroup ${S_0}$ such that
the partial order $\G/{S_0}$ will be a linear order
with a smallest element where  each element has
 a direct successor and
(except for the smallest one) a direct
 predecessor; moreover, all  countable sets will be bounded, and  
 all intervals in this order will be either finite or uncountable.  

The set of growth semigroups (or equivalently, the interval of clones)
above  
${S_0}$ will be the Dedekind completion of this order; it will be a 
linear order with no coatom.

\subsection*{Replacing growth semigroups by filters}

It turns out to be convenient to concentrate on growth functions of a
certain kind,  the functions $h_A$ defined below.

 \begin{Definition}\label{star} \herea16.  \begin{enumerate} \item For
 any infinite set $A \subseteq \N$, we let $h_A$ be defined by
 $$h_A(n) := \min\{a\in A :  a>n \}.$$ 
We let $S_A$ be the growth semigroup generated by~$h_A$.
 \item \label{starstar} 
 For any infinite sets $A \subseteq \N$, and for any  
 growth functions $f,g$ we
write $ f\le_{A} g $ for $f \le_{S_A} g $, or more explicitly: 
$$f\le_{A} g \ :\Leftrightarrow \ 
\exists k: \,\, f \le \xx (\max(g,h_{A}))^k $$ 
 (equivalently: there is a $k$ such that
  $f \le \xx (g\circ h_{A})^k $).  
 \\ 
 Note that there is no point in defining  $\le ^*$, since  
  the statement  
 $$ 
  \exists k: \,\, f \le^* \xx (\max(g,h_{A}))^{k} $$ 
 implies $  f \le \xx (\max(g,h_{A}))^{k'} $ for some
  large enough~$k'$.   
 \item If $U$ is a  filter of  infinite subsets of~$\N$, we let  
 $f \le_U g$ iff there is  $ A\in U$ 
  with $f \le_{A}  g$.  
  \item If $U$ is as above, we let $\G_U:= \{g \in \G:  g
  \le_U \id \}$.   
 \end{enumerate} 
 \end{Definition}

 \begin{Fact} \label{gea.fact}
 \herea18.  
 \ 
 \begin{enumerate} 
 \item  $h_A$ is a growth function.  Conversely,  for any (growth) 
  function $g:\N\to \N$ 
 there is a set $A$ such that $g\le \xx h_A^2 $.  
 \item If $A \subseteq B$, then $h_B \le h_A$, so 
 	 $f \le_B g$ implies  
  $f \le_A g$.  Hence $h_{A\cap B} \ge \max(h_A, h_B)$ (if
 	 $A\cap B$ is infinite). 
\item If $A=^* B$, then the relations $\le_A$ and $\le_B$ coincide. 
 \end{enumerate} 
Therefore: {if} $U$ is a filter of infinite sets, then 
 $\G_U$ is the smallest  growth semigroup containing 
$\{h_A: A\in U\}$. Moreover: the relations 
$$ f \le_{\G_U } g$$
(defined in~\ref{leqS}), and 
$$ f \le_U g $$ 
(defined in \ref{star}(3)) are equivalent. 
 \end{Fact} 
  
  {\sl 
 \subsection*{Outline of the proof of theorem~\ref{thm} }

 The relations $\le_A$ and $\le_U$ are partial
  quasiorders (i.e., reflexive and transitive relations) on~$\G $ .

 Factoring out by the relation  
 $$f\sim_U g  \quad \Leftrightarrow
 \quad      
    f \le_U g \   \&  \ g\le_U f$$ 
 will thus give a partial order~$L_U$.  
  
 We will find a filter $U$ (in fact: $U$ will be an ultrafilter) such
 that   
 the relation $\le_U$ is a {\em linear} quasiorder on~$ \G $, so  
 $L_U:= \G /{\sim_U}$ will 
be linearly ordered. The smallest element in this 
 order is the equivalence class of the function~$\id$, i.e.,  
 the set $\G_U $.  
  
 (It is not  
 necessary  to have this interval linearly ordered in order to obtain 
  a clone $\C$ with $[\C,\O]$ not dually atomic, but it will make  
 the proof more transparent.)  
  
 How do we construct~$U$?   For each pair of growth  
 functions $f,g$ we want  to  
 have a set $A\in U$ witnessing $f\le_A g$ or~$g\le_A f$.    
 This requirement tends to add very ``thin'' sets to~$U$ (since we
 have to   
 add a set $A$ such that some iterate of~$\max(g, h_A)$ dominates $f$) 
  
 Achieving this goal alone (under CH) would be a trivial exercise: Let 
 $(f_i,g_i: i\in \omega_1)$ be an enumeration of all pairs of growth 
 functions, then we can define a $\subseteq^*$-decreasing sequence 
 $(A_i:i\in \omega_1)$ of infinite subsets of~$\N$, making $A_i$ so thin 
 that $f_i \le_{A_i} g_i$ or the converse holds.     We then let  
 $U$ be the filter generated by the sets~$A_i$.  
  
 However, such a naive construction will make the sets $A_i$ so thin  
 that the family 
 $\{h_A: A\in U\}$ might  dominate all unary functions, which would result 
 in~$\G_U = \O$.   
  
 So we will modify this naive construction.  Together with the sets $A_i$ we 
 will construct objects $s_i$ which control how fast the sets $A_i$ may  
 be thinned out.    These auxiliary objects $s_i$ will help to give
a more explicit description of the relation~$ \le_{A_i}$ (see \ref{leq}
 and \ref{not}). 
  }

 \section{Trees and creatures}

 \begin{Definition}\label{def20} 
 \herea20. 
   For any (finite) partial order 
   $(T,{\trianglelefteq}) $ we let  $\leaf(T)$ be the set of maximal 
   elements of~$T$, and $\int(T) := T \setminus \leaf(T)$.  
  
 We let $x\vartriangleleft y $ iff $x \trianglelefteq y$ and $x\not= y$.

  For $\eta\in T$ we let $\succ_T(\eta)$ be the set of direct successors of 
   $\eta$ in~$T$, i.e., $$\succ_T(\eta) = \{ \nu\in T: \eta 
   \vartriangleleft \nu, \lnot \exists \nu'\, (\eta\vartriangleleft \nu' 
   \vartriangleleft \nu)\}$$ 
 (so $\succ_T(\eta)=\emptyset $ iff $\eta\in \ext(T)$)  
  
 We say that $T$ is a tree if~$T$ has a least element, called $\wurzel(T)$ 
 and for every 
  $\nu\in T$ the set $\{\eta: \eta \trianglelefteq\nu\}$ is linearly 
  ordered by~$\trianglelefteq$.  

Elements of a tree $T$ will be called ``nodes'', elements of~$\int(T)$
are ``internal nodes''. 
  
 Elements of~$\ext(T) $  are called ``leaves'' or ``external nodes''. 
  
 For any $\eta\in T$ we let  
 $$ \yy T^\eta := \{ \nu: \eta \trianglelefteq \nu \}$$ 
 (with the induced order).   If $T$ is a tree, then also $\yy T ^ \eta 
 $ is a tree, with root~$\eta$.   
 \end{Definition}

 \begin{Definition} 
 \herea21. 
 Let $(T,{\trianglelefteq})$ be a tree.   
  A {\em branch} of~$T$ is a maximal linearly ordered subset of~$( 
  T,{\trianglelefteq})$.  In other words, a branch is a set of the 
  form $\{\eta: \eta \trianglelefteq \nu\}$ for some $\nu\in \leaf(T)$.   
(We may occasionally identify 
 a  node $\nu\in  \leaf(T)$ with the branch
 $\{\eta: \eta \trianglelefteq \nu\}$.)   
  
 A {\em front} is a subset of~$T$ which meets each branch exactly once. 
 For example, $\leaf(T)$ is a front, and the singleton $\{\wurzel(T)\}$ 
  is also a front.  If $F \subseteq \int(T)$ is a 
 front, then also $\bigcup_{\eta\in F} \succ_T(\eta)$ is a front. 
 \end{Definition}

Note that for any tree $(T, {\trianglelefteq})$, any subset $S \subseteq T$
(with the induced order) will again be a tree. The following
definition singles out some of those subsets.

 \begin{Definition} \label{subtree}
 \herea22. If $(S, {\trianglelefteq})$ 
    is a tree, $T \subseteq S$, we call  
 $T$ a {\em subtree} of~$S$ (``$T\le S$'')
 iff $T$ contains the root of~$S$ and: 
 $\forall\eta\in T\cap \int(S): \emptyset\not= 
 \succ_{T}(\eta)\subseteq \succ_S(\eta)$.  
  
 $(T,{\trianglelefteq})$ will again be a tree, and we have $\ext(T) 
 \subseteq \ext(S)$, $\int(T) \subseteq \int(S)$.   
 \end{Definition}

Below we will need the following version of Ramsey's theorem. 
  
 \begin{Fact}
 \herea25. 
For every natural number  
 $n$ there a natural number $k=\ramsey(n)$ such that:  
 \begin{quote} 
 Whenever $f:[C]^ 2 \to \{0,1\}$, with  $|C|\ge k$, \\
there is 
 a subset $A\subseteq C$ with $|A|\ge n$ which is 
 homogeneous for~$f$, i.e., such that $f\on [A]^ 2 $ is constant. 
 \end{quote} 

Here, $[C]^2 := \{\{x,y\}: x,y\in C, x\not=y\}$ is the set of all
unordered pairs from~$C$.

 It is well known that one can choose $\ramsey(n) = 4^n$ or even smaller, but 
 we do not need any good bounds on~$\ramsey(n)$, its  mere existence is 
 sufficient.   For the rest of the paper we fix such a function~$\ramsey$. 

Note that $\ramsey$ will satisfy $\ramsey(k) \ge 2k$ for all $k\ge 3$.  
We let $\ramsey^{-1}(k) = \max\{n:  \ramsey(n) \le k \}$. 
 \end{Fact}

 \begin{Definition}\label{def28} 
 \herea28. 
  A {\em \bf creature} is a finite tree  $ \T $ 
where 
 \begin{itemize} 
\item 
 $\ext(T) \subseteq \N$, $\int(T) \cap \N = \emptyset$  
 \item Whenever $\eta_1,\eta_2$ are $\trianglelefteq$-incomparable, 
   then either $\max \ext \yy T^{\eta_1} < \min \ext \yy T^{\eta_2} 
   $,  or  
  $\max \ext \yy T^{\eta_2} < \min \ext \yy T^{\eta_1} $.  
 \end{itemize} 
  
 $\T$ is called improper {if} $T$ consists of the single (external)
 node~$\wurzel(T)$.    Otherwise (i.e., {if} $\wurzel(T)\in \int(T)$),
 $\T$  is called   
 a proper creature.

 We will write $\max[\eta]$ and $\min[\eta]$ for $\max \ext\yy T^\eta$  
 and  $\min\ext \yy T^\eta$, 
{if} the underlying tree $T$ is clear from the context. In particular,
 if $\eta\in \ext(T)$ then $\yy T^\eta = \{\eta\}$, so $\max[\eta]=
 \min [\eta]=\eta$. 
  
 If $\eta_1,\eta_2$ are $\trianglelefteq$-incomparable, we will write
 $\eta_1 <_T \eta_2$ or just $\eta_1 < \eta_2$ to abbreviate $\max
 [{\eta_1}] < \min [{\eta_2}]$.  $\eta_1\le \eta_2$ iff $\eta_1 < \eta_2$ or $\eta_1=\eta_2$.  This is a partial order; the order
 $\trianglelefteq$ can be viewed as pointing from the bottom to the top, 
whereas this order $\le$ can be viewed as pointing from left to right. 
  We call $\le$ the lexicographic order of nodes, since clearly 
\begin{quote}
Whenever $\eta_1 < \eta_2$, and   $\eta_1 \trianglelefteq \nu_1$, 
  $\eta_2 \trianglelefteq \nu_2$, \\
then also $\nu_1 < \nu_2$. 
\end{quote}
Note that on $\ext(T)$ this lexicographic order and the usual order of $\N$ agree.

  If $\T$ is a proper creature, then  
 we let $\|\T\| := \min\{ |\succ_T(\nu)|: 
 \nu\in \int(T)\} $.  (If $\T$ is improper, $\|\T\|$ is undefined.) 
 \end{Definition} 

\begin{Remark}\label{470a}
  (Readers not familiar with creature forcing are advised to skip this
 remark.) 
The reader who is familiar with \cite{RoSh:470} will be disappointed 
to see that norms, an essential ingredient in \cite{RoSh:470}'s creatures,
 are not mentioned here at all.  It will become clear that for 
our purposes there is in fact a natural norm, namely 
${\rm nor}_T(\eta):= \max\{k: \xx\ramsey^k (4) \le |\succ_T(\eta)| \}$. 
See also \ref{470b}. 
\end{Remark}

 \begin{Notation} 
 \herea33. 
 \label{coloring} 
Let $E$ be a finite  set $E$ of ``colors'' (often with only 2 elements).  
  
 We will consider three kinds of ``coloring functions'' on creatures~$\S$. 
 \begin{enumerate} 
 \item (Partial)  functions $c: S \to E$.  These are called {\em unary 
  node colorings}.  

 \item (Partial) functions $c : \bigcup_{\eta\in \int(S)}
      [\succ_S(\eta)]^2 \to E$.
  These are called {\em binary node colorings}.   We may write such
 functions more sloppily as $c:[S]^2 \to E$, with the understanding
 that we will only be interested in the values~$c(\{\nu_1,\nu_2\})$, where  
  $\nu_1$ and $\nu_2$ have a common direct $\vartriangleleft$-predecessor.  
 \item  (Partial) functions $c : \ext(S)\to E$, or equivalently, functions $c$
 whose domain is contained in the set of branches of~$S$. 
 These are called {\em unary branch
 colorings}.   
 \end{enumerate} 
 If $c$ is a coloring of one of these types, and $T\le S$, we say that $T$ 
  is {\em $c$-homogeneous} (or  {\em homogeneous for $c$}), iff: 
 \begin{enumerate} 
 \item In the first case:  $c\on \succ_T(\eta)$ is constant,  
 for all $\eta\in \int( T)$.   
 \item In the second case: $c\on [\succ_T(\eta)]^2$ is constant,  
 for all $\eta\in \int(T)$.   
 \item  In the third case: $c\on \ext(T)$ is constant.  
 \end{enumerate} 
 \end{Notation}

 \begin{Lemma} 
 \herea36. 
 \label{homo} 
Let $n\ge 3$, 
and  let $\S$ be a creature with $\|\S\|\ge \ramsey(n)$, let~$E:=\{0,1\}$, and  
 let $c:S\to E $, or  $c:[S]^2 \to E$, or $c:\ext (S)\to E$  
  be a coloring (with  2  colors) of one of the types mentioned above.  
  
 Then there is a creature $\T\le \S$, $\|T\| \ge n $,  which  
 is homogeneous for~$c$.  
  
 Similarly, {if} $c$ is a coloring with at most $2^k$ colors, and $\|S\|\ge
\xx \ramsey^k (n)$, 
then we can find a 
$c$-homogeneous creature $\T\le \S$ with $\|T \| \ge n $.  
 \end{Lemma}

 \begin{proof} 
 We will prove this only for the case of 2 colors.  For $2^k$   colors, apply 
 the argument for  2 colors $k$ times.

 First and second case (unary and binary node colorings): We will 
 define $T$ by (``upward'') 
  induction, starting with $\wurzel(T):= \wurzel(S)$. 
 For any $\eta\in T$ we let $\succ_T(\eta) $ 
 be a $c$-homogeneous subset of~$\succ_S(\eta)$ of size $\ge n$. 
(Such a large homogeneous set exists since $\succ_S(\eta)$ has 
at least $\ramsey(n)$ elements.)   This defines a creature $T$ with
 $\|T\|\ge n$, and it is clear that $T$ is $c$-homogeneous.

 Third case: Let $c:\ext(S)\to \{0,1\}$.  We will define a unary node 
 coloring $c'$ as follows, by (``downward'') induction starting at the 
 leaves: 
 \begin{itemize} 
 \item  If $\eta\in \ext(S)$, then~$c'(\eta)=c(\eta)$, 
 \item If $\eta\in \int(S)$, and $c'\on \succ_S(\eta)$ is already defined, 
 then we choose $c'(\eta) \in \{0,1\}$ such that 
$$ |\{ \nu\in \succ(\eta):  c'(\nu) = c'(\eta) \}| \ge n . $$
This is possible as $|\succ(\eta)| \ge \ramsey(n) \ge 2n$. 
 \end{itemize}

 Now let $\ell_0:= c'({\rm root}(S))$, and define $\T\le \S$ 
 by requiring ${\rm root}(T)={\rm root}(S)$, and  
 $$\forall \eta\in T: \  \succ_T(\eta)  = \{\nu\in \succ_S(\eta):
 c'(\nu)=\ell_0\}$$   
 Then $\|T\|\ge n $, and $c'$ is constant on~$T$ with constant
 value~$\ell_0$, so also $c$ is constant on~$\ext(T)$.

 \end{proof}

 \section{Zoos}

 \begin{Definition}\herea36a. \label{zoodef}
   A {\em zoo} is a sequence  
 $  s = (\S_0, \S_1, \ldots ) $ of proper creatures (see~\ref{def28}) 
 such that  
 \begin{itemize} 
 \item All $S_n$ are pairwise disjoint 
 \item For all $n$, $\max \ext( S_n)  < \min \ext (S_{n+1})$  (recall that 
 $\ext(S_n)\subseteq \N$) 
 \item The sequence $(\|\S_n\|: n\in \N)$ diverges to~$\infty$, and
 $\forall n: \|S_n\|\ge 4$. 
 \end{itemize} 
 We define $\ext(s) := \bigcup \ext(S_n)$, similarly~$\int(s)$.  
  
 We similarly transfer other notations from creatures to zoos, e.g.,  
 for $\eta\in S_n$ we may write $\succ_s(\eta)$ for~$\succ_{S_n}(\eta)$,
$\yy s^\eta$ for $\yy S_n^\eta $, 
 etc. 
 For $\eta\in S_n$, we may  write $\max[\eta]$ or $\max_s[\eta]$ for 
  $\max(\ext(\yy S_n  ^\eta   ))$.  
 We sometimes identify $s$ with   $\int(s)\cup \ext(s)$, i.e.\  
 for  $s=(\S_0, \S_1, \ldots)$  we 
  write $\eta\in s$ instead of~$\eta\in \bigcup_n S_n$.  
 \end{Definition}

\begin{Remark}\label{470b}
(Again this remark is only for the benefit of readers familiar with
\cite{RoSh:470}.)
With the definition ${\rm nor}(A) =
 \max\{k: \xx\ramsey^k (4)\le |A| \}$, it is clear that 
$\lim_{n\to \infty} {\rm nor}(A_n) = \infty$ is equivalent to just 
$\lim_{n\to \infty} |A_n| = \infty$.   This equivalence  allows us to 
omit the definition of norms for our creatures. 
\end{Remark}

 \begin{Definition}\herea36aa. \label{defstar}
 Let $s = (\S_0,\ldots)$ and $t= (\T_0,\ldots)$ be  zoos.  We say $s =^* 
 t$  {if} there are $n_0$ and $k_0$ such that  
 $(\S_{n_0}, \S_{n_0+1}, \ldots) = (\T_{k_0}, \T_{k_0+1}, \ldots)$.  
 \end{Definition}

 \begin{Definition} 
 \herea42. 
   Let  
 $  t = (\T_0, \T_1, \ldots ) $  be a zoo.  
  
 A {\em branch} of~$t$ is a set which is a branch in one of the
 trees~$T_n$.   
  
   A {\em front} is a  
 set $F \subseteq t$ (i.e., $F \subseteq \bigcup_n T_n$)
 such that each of the sets $F\cap T_n$ 
 is a front in~$T_n$.  Equivalently:  $F$ meets every branch of $t$ in
 exactly  one node.

 A {\em *front} (or: {\em ``almost front''\/}) 
 is a set $F \subseteq \bigcup_n T_n$ such that for
 some $t' =^* t$, $F\cap t'$ is a front in~$t'$.  Equivalently, $F$ is
 a *front iff almost all (=all except finitely many) branches
of $t$  meet $F$ in exactly one node.
  
 If $F$ and $G$ are *fronts in~$t$, we write $F \approx^* G$ iff 
 there is some $t'=^*  t $ such that $F\cap t' = G\cap t'$.
 (Equivalently:   
 $F \approx^* G$ iff the symmetric difference of~$F$ and $G$ is finite.)
   
 \end{Definition}

  \begin{Definition}\label{def43} 
 \herea43. 
 Let $s=(\S_0,\S_1,\ldots) $ be a zoo. For $\eta,\nu\in \bigcup_n S_n$
 we define  
 $$ \eta < \nu \ \Leftrightarrow \ \max[\eta] < \min [\nu], $$ 
 and we let $\eta \lessdot \nu $ iff $\eta$ and $\nu$ are on adjacent
 branches, i.e.,   
 \begin{quote} 
  $\eta < \nu$ and for all $k\in \ext( s)$: $k\le \max[\eta]$ or $k\ge
  \min[\nu]$.   
 \end{quote} 
  \end{Definition}

The following fact is not needed, but can be useful to visualize
fronts.

 \begin{Fact} 
 \herea45. 
 Let $s=(\S_0,\S_1,\ldots) $ be a zoo. A set $F \subseteq  
 \bigcup_n S_n$ is a front in~$s$ iff $F$ can be enumerated as  
  $F=\{\eta_0,\eta_1,\ldots\,\}$ where $\min[\eta_0]=\min \ext S_ 0$, and 
 for all $n$ we have $\eta_n \lessdot \eta_{n+1}$.  
 \end{Fact}

Below we will several times have to generate new (``improved'') zoos
from old ones.  The following examples are special cases of the
general definition below.

\begin{Construction}\label{examples}
Let $s = (S_0, S_1, \ldots)$ be a zoo. 
 \begin{description} 
 \item [drop] Let $n_0 < n_1 < \cdots $ be an infinite increasing
 sequence.   Let $t = (T_0,T_1, \ldots)$ be defined by $T_k =
 S_{n_k}$. Then $t$ is a  zoo.  \\
 We say that $t$ is obtained from
 $s$ by dropping creatures (namely, by dropping the creatures $S_i$
 with $i\in \N\setminus \{n_0,n_1,\ldots\}$). 
 \item[shrink]  For each $n$ let $T_n\le S_n$ (see  \ref{subtree}). Then 
 $t = (T_0,T_1, \ldots)$ is a zoo, {\em provided that $\lim_{n\to
 \infty }  \|T_n\| = \infty $ and $\forall n: \|T_n\|\ge 4$.} \\
We say that $t$ is obtained by
 shrinking~$s$. 

 \item[glue] Partition $\N$ into intervals:  
  $\N =\bigcup_{k=0}^\infty [n_k, n_{k+1})$ with  
$0=n_0 < n_1 < n_2 < \cdots $ and 
 $$\lim_{k\to \infty} (n_{k+1}-n_k) =  \infty, \qquad \forall n  :n_{k+1}-n_k>4.$$ 
 For each~$k$, let $T_k := \{\tau_k\} \cup \bigcup_{i=n_k} ^ 
 {n_{k+1}-1} S_i$, where $\tau_k = \wurzel(T_k)$ will be a new element. 
 Make $T_k$ into a creature 
 by requiring:  For all $i\in [n_k,n_{k+1})$, all $\eta\in \S_i$:   
 $\yy\S_i^ \eta = \yy \T_k^\eta$.       Again $t$ will be a zoo.
 \end{description} 
\end{Construction}
Sometimes we have to compose the steps described above, but
 at one point (see \ref{limit}) we will need a more complicated gluing 
 process, as described in the following definitions.

 \begin{Definition}\herea48a. 
 A gluing recipe is a sequence $r=({R}_0, {R}_1, \ldots)$ of (not necessarily 
 proper) creatures  
 satisfying $\forall k: \max \ext( R_k) < \min \ext(R_{k+1})$, and  
 moreover:     If the set $A:= \{k:   \mbox{${R}_k$ proper}\}$ is infinite,  
 then $\lim_{k\in A} \|{R}_k\| = \infty$.  
  
 That is,  a gluing recipe looks like a zoo, except that  
 we allow all or some of the creatures to be just single natural numbers.  
 \end{Definition}

Figure~1 shows a zoo $s$ together with a  gluing recipe~$r$.

 \begin{Definition}\herea49. \label{defle}
 Let $s = (\S_0,\ldots) $ and $t= (\T_0,\ldots) $  be  zoos, and let  
 $r  = ({R}_0, {R}_1,\ldots)$ be a gluing recipe.  
  
 We say that ``$t\le s$ via $r$''  
  iff there are creatures $\S_k'\le \S_k$ with  
 $$\bigcup_n \int(R_n) \ \  \cap\  \  \bigcup_{\kk } \int( S_k')  =
 \emptyset$$   
 such that $t$ is obtained from $r$ by replacing each $k\in \ext(r)$ 
 by~$S_k'$, i.e. 
 $$ T_n = ( R_n \setminus \ext(R_n) ) \ \ \cup\ \  \bigcup_{\kkn} S_k' ,$$ 
and $\yy T_n ^\eta = \yy {{S'_k}}^\eta$ for $\eta\in S_k$, $k\in
\ext(R_n)$, and $\yy T_n^\eta =\int(\yy R_n^\eta) \cup \bigcup _{k\in
\ext(\yy R_n^\eta) } S_k' $ for $\eta\in \int(R_n)$. 

Figure~2 shows the zoo $t$ obtained from the zoo $s$ via the gluing
recipe $r$ from~figure~1.   

 We say ``$t\le s$'' iff there is a gluing recipe $r$ such that  
 ``$t \le s$ via $r$'' holds.  
  
We say  ``$t\le^* s$'' iff there is some $t'=^* t $ such that 
$t' \le s$. 
  
 We leave it to the reader to check that the relation $t\le s$  and
$t\le^* s $ are  indeed  transitive  (and reflexive).   Also,
$\le$ is antisymmetric, and  
$$ s \le^* t \ \ \& \ \ t \le^* s \quad  \Rightarrow \quad s=^* t . $$   
 \end{Definition} 
  
\begin{figure}
\setlength{\unitlength}{0.00035in}
\begingroup\makeatletter\ifx\SetFigFont\undefined%
\gdef\SetFigFont#1#2#3#4#5{%
  \reset@font\fontsize{#1}{#2pt}%
  \fontfamily{#3}\fontseries{#4}\fontshape{#5}%
  \selectfont}%
\fi\endgroup%
{\newcommand{\dashlinestretch}{30}
\begin{picture}(9099,4614)(0,-10)
\put(7062,2512){\ellipse{150}{150}}
\put(537,2487){\ellipse{150}{150}}
\put(3387,2562){\ellipse{150}{150}}
\put(4587,2562){\ellipse{150}{150}}
\path(12,4587)(462,3462)(912,4587)(12,4587)
\path(1587,4587)(2037,3462)(2487,4587)(1587,4587)
\path(3087,4587)(3537,3462)(3987,4587)(3087,4587)
\path(4362,4587)(4812,3462)(5262,4587)(4362,4587)
\path(5637,4587)(6087,3462)(6537,4587)(5637,4587)
\path(6837,4587)(7287,3462)(7737,4587)(6837,4587)
\path(8187,4587)(8637,3462)(9087,4587)(8187,4587)
\path(537,2487)(2187,12)(3687,1812)(3387,2562)
\path(3687,1812)(4587,2562)
\put(307,4105){\makebox(0,0)[lb]{$S_0$}}
\put(3357,4105){\makebox(0,0)[lb]{$S_2$}}
\put(7107,4105){\makebox(0,0)[lb]{$S_5$}}
\put(8457,4105){\makebox(0,0)[lb]{$S_6$}}
\put(9700,4105){\makebox(0,0)[lb]{$\cdots$}}
\put(1857,4105){\makebox(0,0)[lb]{$S_1$}}
\put(4662,4105){\makebox(0,0)[lb]{$S_3$}} 
\put(5907,4105){\makebox(0,0)[lb]{$S_4$}}
\put(687,2487){\makebox(0,0)[lb]{$0$}}
\put(3537,2562){\makebox(0,0)[lb]{$2$}}
\put(4737,2562){\makebox(0,0)[lb]{$3$}}
\put(6937,1900){\makebox(0,0)[lb]{$R_1$}}
\put(1887,1900){\makebox(0,0)[lb]{$R_0$}}
\put(7212,2512){\makebox(0,0)[lb]{$5$}}
\put(8700,1700){\makebox(0,0)[lb]{$\cdots$}}  %
\end{picture}
}
\begin{center}Figure 1.
\end{center}
 
  \end{figure}

\begin{figure}

\setlength{\unitlength}{0.00083333in}
\begingroup\makeatletter\ifx\SetFigFont\undefined%
\gdef\SetFigFont#1#2#3#4#5{%
  \reset@font\fontsize{#1}{#2pt}%
  \fontfamily{#3}\fontseries{#4}\fontshape{#5}%
  \selectfont}%
\fi\endgroup%
{\newcommand{\dashlinestretch}{30}
\begin{picture}(3999,2052)(0,-10)
\path(312,1275)(12,2025)(612,2025)
	(312,1275)(687,375)(1512,900)
	(1287,1350)(1062,2025)(1512,2025)(1287,1350)
\path(1512,900)(2187,1350)(1962,1950)
	(2412,1950)(2187,1350)
\path(3187,600)(2737,1950)(3487,1950)(3187,600)
\put(237,1725){\makebox(0,0)[lb]{$S_0'$}}
\put(1212,1725){\makebox(0,0)[lb]{$S_2'$}}
\put(2112,1650){\makebox(0,0)[lb]{$S_3'$}}
\put(3012,1425){\makebox(0,0)[lb]{$S_5'$}}
\put(612,50){\makebox(0,0)[lb]{$T_0$}}
\put(3112,50){\makebox(0,0)[lb]{$T_1$}}
\put(3912,1425){\makebox(0,0)[lb]{$\cdots$}}
\end{picture}
}
\begin{center}Figure 2.
\end{center}
\end{figure}

For visualizing creatures, and also for avoiding notational complications, 
it is often useful to replace the relation $\le$ by the following relation~$\leqq$:   
  
 \begin{Convention}\label{defqq}
 We will write  $t\leqq s$ {if} $t\le s$ via some gluing recipe~$r$, and
 in addition to   
 $$\bigcup_n \int(R_n )\ \   \cap\ \  \bigcup_{\kk} \int( S_k')  =
 \emptyset$$   
 we have  moreover  
 $$\bigcup_n \int(R_n) \cap \bigcup_n \int(S_n) = \emptyset, $$ 
  i.e., the internal nodes from $s$ that  
 were omitted (either in the passage from $S_n$ to~$S_n'$, or because they 
 are in some $S_k$ with $k\notin \ext(r)$) 
  will not be recycled as nodes from~$r$.   
 Note: 
 \begin{enumerate} 
 \item  The relation $\leqq$ is not  transitive.  
 \item  {\em However:} For any $t\le s$ we can find (by renaming
  internal nodes  
  of~$t$) a zoo $t' \leqq s$, $t' = (\T_0', \ldots)$ 
  such that $t'$ is isomorphic to~$t$, i.e.,  
 there is a bijection between $\bigcup _n T_n'$ and $\bigcup_n T_n$ 
 which is the identity on leaves, and preserves the relations
  $\trianglelefteq$   and the norms.  
 \item {\em Moreover:} 
 Let $(s_i: i \in I)$ be a family of zoos, and $\forall i: t\leqq s_i$. 
 Assume $t_1 \le t$ via~$r_1$.  Then there is a gluing recipe $r_1'$ which 
 is isomorphic to~$r_1$ 
 and a condition $t_1'$ isomorphic to~$t_1$, with $t_1\leqq t $ via~$r_1'$, and 
 also satisfying $\forall i: t_1'\leqq s_i$. 
 \end{enumerate} 
 Thus, whenever we consider conditions $t\le s$ we will usually assume
 without loss of generality that we have even $t \leqq s$.  This
 guarantees that any $\eta \in t \cap s$ will appear in~$t$ ``in the
 same place'' as in~$s$, e.g., $\succ_t(\eta) \subseteq \succ_s(\eta
 )$.

We write $t \leqq^* s $ iff there is $t' =^* t$ with $t' \leqq s$. 

 \end{Convention}

\begin{Definition} \label{color.zoo}
Let $E$ be a finite set, and
let $t$ be a zoo. 
A coloring (of~$s$, with colors in $E$)
 is  a partial map $c: s \to E$, or $c:[s]^2 \to E$, or 
$c:\ext(s)\to E$, such that each map $c\on S_n$, or $c\on [S_n]^2$, 
or $c\on\ext(S_n)$ is a coloring as in \ref{coloring}.  Again we call
$c$ a unary node coloring, a binary node coloring, or a branch
coloring, respectively. 

Let $t\leqq s$. 

We say that $t\leqq s$ is $c$-homogeneous if: 
\begin{itemize} 
 \item In the first case:  for all~$n$, $c\on \succ_{T_n}(\eta)$ is constant,  
 for all $\eta\in \int( T_n)\cap s$. 
 \item In the second case: for all~$n$, 
$c\on [\succ_{T_n}(\eta)]^2$ is constant,  
 for all $\eta\in \int(T_n)\cap s$.   
 \item  In the third case: not only is each $c\on \ext(T_n)$ is
 constant, but all constant values are the same, i.e.: $c\on \ext(t)$
 is constant. 
\end{itemize}

We say that $t \leqq^* s $ is almost $c$-homogeneous iff there is 
$t' =^* t $, such that $t'$  is $c$-homogeneous. 

  \end{Definition}

\begin{Fact} \label{homo.persists}
 If $t$ is $c$-homogeneous, then any $t'\leqq t$ is also $c$-homogeneous.  

  If $t$ is almost $c$-homogeneous, 
then any $t'\leqq^* t$ is also almost $c$-homogeneous.  
\end{Fact}

 \begin{Lemma} 
 \herea80. 
 \label{homop} 
 Let $s=(\S_0,\ldots) $ be a zoo, and let $E$  be a finite set 
   Let $c:\bigcup_n S_n \to E $,  
 or  
         $c:\bigcup_n [S_n]^2 \to E $, 
 or 
  $c:\bigcup_n \ext(S_n) \to E $ be a coloring of $s$ with 
colors from~$E$.   
  
 Then there is a zoo $t\leqq s$ which is homogeneous for~$c$. 

 Moreover: $t$ can be obtained from $s$ by combining the steps
``shrinking'' and ``dropping'' (i.e., with gluing recipes which
contain only improper creatures).
 \end{Lemma}
  
 \begin{proof} 
We show this only for the case~$|E|=2$. (For larger~$E$, repeat the 
proof $\lceil\log_2 |E|\rceil $ many times, or use the unproved
assertion from lemma~\ref{homo}.)

We may assume $\|S_n\|\ge \ramsey(4)$ for all~$n$.   By Ramsey's theorem
we can find a sequence $(\ell_n:n\in \N) $ [namely: $\ell_n:= 
\ramsey^{-1} ( \|S_n\|)$] which diverges
to infinity and satisfies $\forall n: \|S_n\| \ge \ramsey(\ell_n)$.

 Apply lemma~\ref{homo} to each $S_n$ separately to get  creatures
 $T_n\le S_n$ which are homogeneous for~$c$ and satisfy $\|T_n\|\ge
 \ell_n$. Thus $t= (T_0,T_1,\ldots)$ is a zoo. 
  
  If $c$ is a coloring of the third kind (a unary branch coloring),
 then it is still possible that the constant values that $c$ takes on
 each creature are different.  One of the constant values appears
 infinitely often, so by dropping creatures from $t$ we obtain
 $t'\leqq t$ which is $c$-homogeneous.
 \end{proof}

It is clear that we can extend this lemma to the slightly more general case of finitely many colorings:  
If  $s=(\S_0,\ldots) $ is a zoo, and~$c_1$, \dots, $c_k$ are
 colorings of~$s$, each with finitely many colors $E_1$, \dots, $E_k$,
respectively,
 then there is a 
zoo $t\leqq s$ which is homogeneous for each~$c_i$. 

If we have countably many colorings then we can in general not find a
zoo which is homogeneous for all of them; however, the following
construction shows that this is almost possible.

 \begin{Lemma}\herea243. 
  \label{limit}  Let $(s_0, s_1,\ldots)$ be a sequence 
 of zoos with $\forall k,n: k< n \Rightarrow  s_{n} \leqq^* s_k$.  
  
 Then there is a zoo $t$ with $\forall n: t\leqq^* s_n$.  
 \end{Lemma} 
  
 \begin{proof} 
  
 Let $s_n = (\S^n_0, \S^n_1, \ldots)$.   Without loss of generality 
 we may assume (omitting finitely many creatures from 
 $s_{n}$ {if} necessary) that we have in fact $s_{n} \leqq s_k$ for  
 all~$k<n$.  
  
 Moreover we may assume, for all~$n$: 
\begin{enumerate}
\item [(a)] 
$\min  \ext S^{n+1}_0   > \max \ext S^n_0$, and  
\item[(b)] $\|\S_0^n\| \ge n+4$. 
  \end{enumerate}
 Now let $t:= (\S^0_0, \S^1_0, \S^2_0, \ldots)$, then we claim 
that $t\leqq^* s_n$
 for all~$n$.  

First, note that $t$ is indeed a zoo, by (a) and (b).

We will show only $t \leqq s_0$; a similar proof establishing 
$(\S^n_0, \S^{n+1}_0, \ldots)\leqq s_n$ is left to the reader. 

Write $S_n$ for~$S_n^0$. 

Since $s_n = (S^n_0, S^n_1, \ldots \,) \leqq s_0$, there is a creature
$R_n$ such that $S^n_0 = \int(R_n) \cup \bigcup_{k\in \ext(R_n)} S_k'$, 
where
$S_k' \le S_k$.  Note $\|R_n\| \ge \| S^n_0\| \ge n$ (if $R_n$ is proper). 

We leave it to the reader to check that $\max \ext R_n < \min \ext R_{n+1}$
(using $\max \ext S ^n_0  < \min \ext S^{n+1} _0$). 

Now the gluing recipe  $ r = (R_0,R_1,\ldots) $ witnesses $t \leqq s_0$. 
 \end{proof}

The following corollary will be in our transfinite
 construction~\ref{CH1}.  See also  
 \ref{finalremark}.

 \begin{Corollary}\herea243a. 
  \label{limit2}  Let $Z$ be a  countable set of zoos which is 
linearly quasiordered by~$\leqq^*$.  Then 
there is a zoo $t$ with $\forall s\in Z: t\leqq^* s$.  

 \end{Corollary}
\begin{proof}  Let $Z= \{z_0, z_1, \ldots\}$.   Choose a 
sequence $s_0\geqq^* s_1 \geqq^* \cdots $ in $Z$ such that $s_n\in Z$,
and $s_n \leqq^* z_n$, 
and apply lemma~\ref{limit}. 
\end{proof}

\begin{Corollary}
 \label{homoc} 
 Let $s=(\S_0,\ldots) $ be a zoo, and  for each $n\in \N$ let $c_n:s\to E_n$ 
or $c_n:[s]^2\to E_n$ or $c_n:\ext(s)\to E_n$ be a coloring with finitely many 
colors. 
  
 Then there is a zoo $t\leqq^* s$ which is almost homogeneous for each~$c_n$.  That is:  
 there is a zoo~$t\leqq^* s$, and for each $n$ there is $t_n=^* t$, $t_n\leqq s$, which is homogeneous  for~$c_n$.   
\end{Corollary}
   
\begin{proof} Apply lemma~\ref{homop} infinitely many times to get a
decreasing sequence $s \geqq s_0 \geqq s_1 \geqq \cdots $, such that each $s_i$ is $c_i$-homogeneous.   By \ref{limit} there is a  $\leqq^*$-lower bound $t$
for this sequence;  by \ref{homo.persists}, $t$ is almost $c_i$-homogeneous for each~$i$. 
  \end{proof}

\section{Gauging growth}

\begin{Motivation}
{\sl 

The comparison ``$f \le g$'' between growth functions is too coarse for
our purposes. For finer comparisons, we will consider the growth
behavior of these functions
``locally''.  For example, {if} $I = [a,b] $ is an interval in~$\N$, and
$f(I)\cap I = \emptyset $ (i.e., $f(a) > b$), then we can say that $f$
grows fast on~$I$, or symbolically: $f\on I > I$. 
Similarly, {if} $\I = \{ [a_1, b_1], \ldots, [a_n ,b_n]\}$
is a set of intervals with $
a_1 < b_1 < \cdots <  a_n < b_n$, and $f(a_1) > b_n$, we can say that
$f$ 
``grows faster than $\I$'', or ``is
stronger than the set $\I$'', symbolically: $f\on \I > \I$. 

This point of view allows us to introduce a dual concept:  {if} 
   again   we have $a_1< b_1 < \cdots < a_n < b_n$, but now
 $f(b_k) < a_{k+1}$ for $k=1,\ldots, n-1$, then we can say that 
$f$ ``grows more slowly than $\I$'', or:
$f$ is ``weaker''  than the set 
$\I$, symbolically: $f\on \I < \I$.

Note that this is indeed a local notion:  If $\I_1$ and $\I_2$ are
sets of intervals, $f$ and $g$ growth functions, then it is quite possible 
that $f\on \I_1 < \I_1$  but $f\on \I_2 > \I_2$, while $g$ satisfy the
converse inequalities. 

How does this help us to compare $f$ and~$g$? 
If we label certain sequences of intervals as `distinguished', we 
introduce a (kind of) ordering relation
  on growth functions:  $f < g$ iff there is a 
distinguished sequence of intervals
$\I$  which is stronger than $f$ but weaker 
than~$g$: 
$f\on \I < \I <  g\on \I$. 
     Our aim in the remaining sections is to show
how we can select such distinguished sequences such that the resulting
ordering relations can be viewed as a linear order. 

Zoos (or rather: nodes in zoos)  are our way of coordinatizing
sequences of intervals.    If $s$ is a  zoo, $\eta\in s$, then 
$\eta$ is associated with the interval $[\min[\eta], \max[\eta]]$, 
and also with the set of intervals  $\{ [\min[\nu], \max[\nu]]: 
\nu\in \succ(\eta)\}$. 
}
\end{Motivation}

 \begin{Definition} \label{sw}
 \herea84. 
  Let $f$ be a growth function,  
    $s= (\S_0,\S_1, \ldots)$ a zoo.  
 \begin{itemize} 
  
 \item[(a)]  {\em $s$ is $f$-strong} iff 
  $$ \forall n: \ f(\max \ext( S_n)) \le \min \ext( S_{n+1}) $$ 
  Not surprisingly, we say that $s$ is {\em almost $f$-strong} if the above
inequality holds for all but finitely many~$n$. 
 \item[(b)]  Let $\eta\in \int(S_n)$.   We call $\eta$ 
 	{\em $f$-strong in~$s$} iff: 
 \begin{quote} 
  For all $\nu_1<\nu_2$ in~$\succ_{S_n}(\eta):$ 
   $f(\max [\nu_1]) \le \min [\nu_2] $ 
 \end{quote} 
 \item[(c)]  
  Let $\eta\in S_n$.   We call $\eta$ 
 	{\em $f$-weak in~$s$} iff: 
 \begin{quote} 
   
  $f(\min [\eta] ) > \max [\eta] $ 
 \end{quote} 
 \end{itemize} 
 \end{Definition} 
  
 \begin{Fact} \label{noboth}
 \herea88. 
\begin{enumerate} 
\item   No $\eta$ can be both  
 weak and strong. 
\item    For $\eta \in \ext(S_n)$, $\max[\eta]=\min[\eta] = 
\eta\in  \N$.
 So (since $f$ is a growth function), every leaf is $f$-weak.
\item If $\eta\trianglelefteq \nu$, and $\eta$ is weak, then $\nu$ 
 cannot be strong.    
\item  If $t\leqq s$, and $\eta\in t\cap s$ is $f$-weak in~$s$,
 then $\eta$ is also $f$-weak in~$t$.   
\item 
Similarly, if  $\eta$ is $f$-strong  in~$s$,
 then $\eta$ is also $f$-strong in~$t$.   
\item Finally, if $s$ is $f$-strong, then any $t\leqq s$ is also $f$-strong.
\end{enumerate}
 \end{Fact}

\begin{Fact}  Let $s$ be a zoo, $f$ a growth function.  Then there is 
$t \leqq s$ which is $f$-strong. 
\end{Fact}
\begin{proof} 
The sequence 
  $(\min(\ext(S_n)): n=0,1,2,\ldots)$ 
  diverges to infinity, so  
 replacing $(\S_0,\S_1,\ldots)$ by a 
 subsequence, {if} necessary, we get $f(\max S_n) < \min 
 (S_{n+1})$ for all~$n$.    
\end{proof}

 \begin{Definition} \label{front}
 Let $s$ be a zoo, $f$ a growth function, $F \subseteq s$ a *front.  
  
 We say that $F$ gauges $f$ (or more precisely: $F$ gauges $f$ in~$s$)   iff 
 \begin{enumerate}   
 \item All $\eta\in F$ are $f$-weak 
 \item Whenever $\eta_1  < \eta_2$ are in~$F$, then $f(\max[\eta_1]) \le 
 \min[\eta_2]$. 
\item $s$ is almost $f$-strong. (This actually follows from (1).)
 \end{enumerate} 
  
 We say that $s$ gauges $f$ iff there is a *front  $F \subseteq s $ 
 which gauges~$f$.    In this case we fix such a *front
(or, if possible, an actual front)  and call it~$
 F(s,f)$.  
 \end{Definition}

 The following lemma is easy but important:  
  
 \begin{CL} \herea100. \label{crucial}
 \begin{enumerate} 
 \item If $F_1, F_2  \subseteq s$ both gauge~$f$, then $F_1 \approx^* F_2 $.  
 \item If $s$ gauges~$f$, and  $t\leqq^* s$, then also $t$  gauges~$f$,  
  and $F(t,f) \approx^* F(s,f)\cap t$.  
 \end{enumerate} 
 \end{CL} 
 \begin{proof}

 (1): Without loss of generality (dropping finitely many creatures {if} 
    necessary) we may assume that $F_1$ and $F_2$ are not only *fronts 
    but actually fronts.  Now assume $F_1\not= F_2$.  So there is 
    (wlog) $\eta\in \int(T_n) \cap ( F_1\setminus F_2)$ such that $F_2\cap \yy 
    T_n^\eta $ is a front in~$\yy T_n^\eta $.  Let $\nu_1 < \nu_2$ be 
    in~$F_2\cap \yy T_n^\eta $.  Since $F_2$ gauges $f$ we must have 
    $f(\max[\nu_1]) \le \min(\nu_2)$, but since $\eta$ is $f$-weak we 
    must also have $f(\min[\eta] ) > \max[\eta]$.  Clearly 
    $\min[\eta]\le \min[\nu_1]$ and $\max[\eta]\ge \max[\nu_2]$, so we 
    have a contradiction. 
  
 (2): Let $F':= F(s,f)\cap t$.  Clearly $F'$ is a *front in~$t$, and
    $F'$ gauges~$f$, so $F' \approx^* F(t,f)$ by (1).
  
 \end{proof} 

The following fact is easy:

 \begin{Fact} 
 \herea126. 
 Let $F$ be a front in~$s$.  Then there is a growth function
 $f$ and $t \leqq s$
such that  
 $t$ gauges  $f$,  and $F(t,f) \approx^* F\cap t$. 
  (In fact we can choose~$t=s$.)
  \end{Fact} 
 \begin{proof} Choose a growth function $f$ satisfying $f(\min[\eta])
 = f(\max[\eta]) = \max[\eta]+1$ for all $\eta \in F$.   We leave the
 details 
 to the reader. 
 \end{proof}

 \begin{Lemma}\label{successor1} \herea102.  Let $f$ be a growth
 function, $s= (\S_0,\S_1, \ldots)$ a zoo which is almost $f$-strong.
 Then there is $t\leqq s$
such that $t$ gauges~$f$. 

In fact, we will find a coloring function $c$ such that any 
$t\leqq s$ which is almost
 $c$-homogeneous will gauge~$f$, and then invoke 
\ref{homop} to show that there is such~$t$.  
 \end{Lemma}

\begin{proof} 

We start with a zoo 
    $s= (\S_0,\S_1, \ldots)$ which is $f$-strong.

   By ignoring  finitely many of the $\S_i$ we may assume $\|\S_n\| 
   \ge\ramsey( \ramsey(4))$ for all~$n$.  
  
  For each  
 $\eta\in \int(s)$ we define a pair coloring $c_\eta$ of~$\succ_s(\eta)$ with  
 three colors  
as 
 follows:   Whenever $\nu<\nu'$ in~$\succ_s(\eta)$, then  
 \begin{itemize} 
  
 \item $c_\eta\{\nu, \nu'\} = {\tt strong}$, {if} $f(\max [\nu]) < 
  \min [{\nu'}] $.  
 \item $c_\eta\{\nu, \nu'\} = {\tt weak}$, {if} $f(\min [\nu]) > 
   \max [{\nu'}] $.  
 \item $c_\eta\{\nu, \nu'\} = {\tt undecided}$, otherwise.  
 \end{itemize} 
  
 Note:  If $\nu_1 < \nu_2 < \nu_3 < \nu_4$ are in~$\succ_S(\eta)$, then  
 at least one of  
 $$ c(\{\nu_1,\nu_4\}) = {\tt strong} \  
 \qquad  \mbox { or } \qquad  
  c(\{\nu_2,\nu_3\}) = {\tt weak }  $$ 
 has to hold, since otherwise we would have  
 $$ f(\max[\nu_1]) \ge \min [\nu_4] \mbox{ \ and \ } f(\min[\nu_2]) \le 
 \max[\nu_3],$$ which together with 
$$
\min[\nu_1] \le \max[\nu_1] <
\min[\nu_2] \le \max[\nu_2] <
\min[\nu_3] \le \max[\nu_3] <
\min[\nu_4] \le \max[\nu_4] $$
yields a contradiction to the fact that $f$ is monotone. 

(Also 
 note: {if} $f(\max(\ext\yy S^{\nu} )) < \min(\ext \yy S^{\nu'} )$, 
and $T\le S$ 
 with $\nu, \nu'\in T$, then also $f(\max(\ext \yy T^{\nu} )) < \min(\ext\yy 
 T^{\nu'} )$.)
  
The family $(c_\eta: \eta\in \int(s))$ defines on $s$ a binary node
coloring~$c$.  Let $t\leqq s$ be $c$-homogeneous, $t= (T_0,T_1,\ldots)$, 
$\|T_n\|\ge 4$ for all~$n$.  

 Since each set $\succ_{T_n}(\eta)$ has
 more  than 3 elements, it is impossible that $c_\eta$ is constantly
 ``{\tt  undecided}''.

 Clearly each $\eta\in t$ is either
  $f$-weak or $f$-strong.

\bigskip
Now we show that any $c$-homogeneous zoo $t \leqq^*$
gauges~$f$. (Note: if $t\leqq s$  via $r$, 	then the $f$-strength 
of $s$ ensures that almost all $\eta\in \int(r)$ will be $f$-strong in
$t$.)

 On every branch $b$ let $\eta_b$ 
  be the $\vartriangleleft$-lowest node which is 
  $f$-weak (recall that all leaves are $f$-weak), 
  and let $F:= F(t,f) := \{ \eta_b:  b \mbox { a branch in~$t$ }\}$.  
  
 Note that {if} $\nu \trianglerighteq \eta_b$, then $\nu$ is $f$-weak 
(by fact~\ref{noboth}), 
while   any $\nu\vartriangleleft \eta_b$ is $f$-strong. 

So $F\cap b =    
 \{\eta_b\}$ for all~$b$.   Hence $F(t,f)$ is a front.
  
 Let $\eta_1 < \eta_2$ be in~$F$, $\eta_1\in T_{n_1}$, $\eta_2\in T_{n_2}$. 
We have to check that $f(\max[\eta_1]) \le \min[\eta_2]$. 

The case $n _1< n_2$ is trivial (since $t$ is $f$-strong). 

So assume~$n_1=n_2 =:n$.  In $(T_n,
 {\trianglelefteq})$ let $\nu$ be the greatest lower 
bound of~$\eta_1 $
 and $\eta_2$.   Then $\nu \vartriangleleft \eta_1  $, so 
$\nu$  is $f$-strong.  Let $ \nu_1 <\nu_2$ in~$\succ(\nu)$, $
 \nu_1\trianglelefteq \eta_1$, $ \nu_2\trianglelefteq \eta_2$.  Clearly
 $f(\max[\eta_1])\le f(\max[\nu_1]) \le \min [\nu_2] \le \min[\eta_2] $
 (where the middle inequality holds because $\nu$ is $f$-strong).
  
 Hence $F$ gauges~$f$.  
  
\end{proof}

{\sl 
We now fix a zoo~$s$; the fronts in~$s$ are naturally partially
ordered by the relation ``is everywhere higher''.  We will show below
that 
the relation $f\le_{\ext(s)} g$ (see~\ref{star}) can be
translated to a ``$F(s,f)$ is higher than $F(s,g)$'', for 
sufficiently  small~$s$. 
}

 \begin{Definition} 
 \herea105.\label{prec}   
 Let $s$ be a zoo, and let $F$ and $G$ be  *fronts 
  in~$s$. 
 We write $F \prec G$   
 (or $F\prec_s G$)  
 iff, for all  
  branches $b$ of~$s$,  
 $b$ meets $F$ $\vartriangleleft$-{\bf above}~$G$, i.e.,  
 letting  
 $b \cap F = \{\eta_{F,b}\}$,  
 $b \cap G = \{\eta_{G,b}\}$, we have: 
  $\eta_{F,b}  \vartriangleright\eta_{G,b}$.  
    (The reason for this apparent reversal of inequalities will become 
    clear in remark~\ref{why} below). 
  
  Recall that $F \approx^* G$ iff 
 there is some $t'=^*  t $ such that $F\cap s = G\cap s$.  

 We write $F \prec_s^* G$ {if} the relation $\eta_{F,b} 
  \vartriangleright\eta_{G,b}$ holds for {\em  almost all} branches $b$ of~$s$ 
  (i.e., for all except finitely many).  Equivalently, $F \prec_s^* G$ 
  iff there is some $s' =^* s$ such that  
 $F \cap s' \, \prec_{s'} \, G\cap s'$.  
  
 Similarly we define  $F \preccurlyeq_s G$ [$F \preccurlyeq_s^* G$,  ] by  
 requiring that  
  for all [except finitely many]  
 branches~$b$, $b$ meets $F$ $\trianglelefteq$-{\bf above}~$G$, i.e.,  
 letting  
 $b \cap F = \{\eta_{F,b}\}$,  
 $b \cap G = \{\eta_{G,b}\}$, we have: 
  $\eta_{F,b}  \trianglerighteq\eta_{G,b}$.

 We write $F + 1 \approx^*_s G$ iff $F \prec_s^* G$ and moreover,  
  for almost all branches $b$ of~$s $, $\eta_{F,b}$ is a direct successor  
 of~$\eta_{G,b}$.   
  
  Similarly we write 
  $F + n \approx^*_s G$
 iff, 
  $F \preccurlyeq_s^* G$ and moreover,  
  for almost all branches $b$ of~$s $, $\eta_{F,b}$ is exactly  $n$
  nodes 
  above~$\eta_{G,b}$.     

The notations $F+n \preccurlyeq_s^* G$  and  
$G \preccurlyeq_s^* F+n$  have the obvious meanings.  In particular, 
$F+1 \preccurlyeq^*_s G$ will be equivalent to 
$F \prec^*_s G$.

 Finally, we write $F + \infty \preccurlyeq_s ^* G$ iff for all~$n$, 
  for almost all branches $b$ of~$s$, $\eta_{F,b}$ is more than $n$ 
  nodes  above~$\eta_{G,b}$.  (Equivalently:  {if} for all~$n$, 
 $F + n \preccurlyeq_s ^* G$.) 
 \end{Definition} 

It is easy to see that [$F+n \approx^* G$ and $F+n \approx^* G'$]
implies $G \approx^* G'$, so this functional notation is
justified.   

However: There are fronts $F$ such that $F+1$ is undefined.  For example
if we let $F$ be the front 
$$ \wurzel(s):= \{ \wurzel(S_0), \wurzel(S_1), \ldots \}$$ 
then  there is no *front $G$ with $F+1 \approx^*_s G$. 


 \begin{Fact}\herea105a. \label{persist}
 If $F \preccurlyeq_s^* G$, and $t\leqq^* s$, then also  
 $(F\cap t) \preccurlyeq_t^* (G\cap t)$.   Similarly for 
  $F + n \preccurlyeq_s^* G$ or    
  $F + \infty \preccurlyeq_s^* G$. 
 \end{Fact}

\begin{Lemma}\label{dichotomy} \herea 105b.
\begin{enumerate}
\item Let $s$ be a zoo, and let $F,G$ be *fronts in $s$. Then there is
a coloring $c$ such that: whenever $t\leqq^* s$
 is almost $c$-homogeneous, then
$F \preccurlyeq_t^* G$ or $G \preccurlyeq_t^* F$. 
\item Let $s$ be a zoo, and let $F\preccurlyeq_s^*G$ be *fronts in $s$.
 Then there is a  coloring  $c$ of $s$ such that:
 whenever $t\leqq^*s$ is almost   $c$-homogeneous, then
 exactly one of the following holds: 
\begin{itemize}
\item  $(F\cap t)  \approx^*_{t} (G\cap t)$ 
\item   $(F\cap t)  + 1 \preccurlyeq^*_{s_j} (G\cap t)$.  
\end{itemize}
\item Let $s$ be a zoo, and let $F\preccurlyeq_s^*G$ be *fronts in $s$.
 Then there are colorings $c_0, c_1, \ldots$ such that:
 whenever $t\leqq s$ is almost $c_k$-homogeneous for all $k$,  then
 exactly one of the following holds: 
\begin{itemize}
\item  $(F\cap t) + n \approx^*_{t} (G\cap t)$ \ for some (unique) $n$
\item   $(F\cap t)  + \infty \preccurlyeq^*_{s_j} (G\cap t)$.  
\end{itemize}
\end{enumerate}
\end{Lemma}

\begin{proof}

(1) For each branch $b$ of $s$ that meets $F$ ($G$, respectively) in a
unique point, let $\{\eta_{F,b}\} = F \cap b$ ($\{\eta_{G,b}\} = G
\cap b$, respectively).

Now color each branch $b$ as follows: 
\begin{itemize} 
\item
$c(b) = \mbox{\tt small}$ if 
                  $\eta_{F,b} \trianglerighteq \eta_{G,b} $
\item 
$c(b) = \mbox{\tt big}$ if 
                  $\eta_{F,b} \vartriangleleft \eta_{G,b} $
\item 
$c(b) = \mbox{\tt unknown}$ if 
                  $\eta_{F,b}$ and/or  $\eta_{G,b}$ is undefined. 
\end{itemize}
Now let $t \leqq^*s$ be almost homogeneous for $c$.  Clearly the color
{\tt unknown} appears only finitely many times as a value of $c\on\ext(t)$.
If $c\on \ext(t)$ is almost constant with value {\tt small}, then
 $F\cap t\preccurlyeq^*G\cap t$, otherwise  $G\cap t\preccurlyeq^* F\cap t$.

\bigskip  
(2) Define $\eta_{F,b} $ and $\eta_{G,b}$ as above.
    Define  a  branch coloring $c$ as follows: 
\begin{itemize} 
\item
$c(b) = \mbox{\tt equal}$ if 
                  $\eta_{F,b}  =  \eta_{G,b} $
\item 
$c(b) = \mbox{\tt bigger}$ if 
                  $\eta_{F,b} \vartriangleright \eta_{G,b} $
\item 
$c(b) = \mbox{\tt unknown}$ if 
                  $\eta_{F,b}$ and/or  $\eta_{G,b}$ is undefined, or if 
                   $\eta_{F,b} \vartriangleleft \eta_{G,b} $
\end{itemize}
Again any almost homogeneous condition can take the value {\tt unknown} at most
finitely many times.
\bigskip

(3) For each $n=0,1,2,\ldots$ define a branch coloring $c_n$ as follows: 
\begin{itemize} 
\item
$c_n(b) = \mbox{\tt small}$ if $\eta_{F,b}$ is at most $n$ nodes
                  $\vartriangleleft$-above $\eta_{G,b} $ 
\item 
$c_n(b) = \mbox{\tt big}$ if 
                  $\eta_{F,b}$ is more than $n$ nodes above $\eta_{G,b} $
\item 
$c_n(b) = \mbox{\tt unknown}$ if 
                  $\eta_{F,b}$ and/or  $\eta_{G,b}$ is undefined, or if 
                   $\eta_{F,b} \vartriangleleft \eta_{G,b} $
\end{itemize}

Now assume that $t \leqq^* s$ is almost $c_n$-homogeneous, for all $n$. 
We distinguish two cases: 

\begin{itemize}
\item[(Case 1)] There is some $n$ such that $c_n\on \ext(t)$ is
constantly {\tt small} (with finitely many exceptions).  Let $\bar n$ be
the smallest $n$ for which this happens, then on almost all branches  $b$ 
of $t$, 
$\eta_{F,b}$ is exactly $\bar n$ nodes above $\eta_{G,b}$, so 
$$ (F \cap t) + \bar n \approx^*_t (G \cap t)$$
\item[(Case 2)] Each $c_n$ is (almost equal to) the constant 
function with value  {\tt big}.  Then we can easily see that 
 $$ (F \cap t) + \infty \preccurlyeq^*_t (G \cap t)$$
\end{itemize}
\end{proof}

 \begin{Definition}\label{precf} 
 Let $f$ and $g$ be growth functions, and assume that $s$ gauges both  
 $f$ and~$g$. (So $F(s,f)$ and  
 $F(s,g)$ are well-defined.)  
  
 We now write $f \approx^*_s g$, $f \preccurlyeq^*_s g$, etc., iff $F(s,f) \approx^*  
  F(s,g)$, $F(s,f) \prec^*_s F(s,g)$, etc, respectively.  

 \end{Definition}

Combining \ref{crucial} and \ref{persist}, we get:  
 If  
 $ f \prec^*_s g $, and $t \leqq^* s$, then also $ f \prec^*_t g $, etc.

 \begin{Remark}\label{why} 
 \herea106.  If $f \le g$, then every $g$-strong node is also $f$-strong,  
 and every $f$-weak node is $g$-weak.  
  
 Hence the front corresponding to~$f$ is {\em $\trianglelefteq$-higher} in the 
 trees than the one for~$g$.   
  
 Thus: $f \le g$ implies $f \preccurlyeq_s g$, whenever $s$ gauges $f$
 and $g$.

 \end{Remark} 

The converse is of course not true, but we will show below that it is
``true modulo $\ext(s)$'':

 \begin{Lemma}\label{leq} 
 \herea121. Let $s$ be a zoo, and let $f$ and $g$ be 
 growth functions, and let  $s$ gauge $f$ and~$g$.  
  
  Assume $f \preccurlyeq_s^* g$.

  Then   $f  \le_{\ext(s)} g  $.  (See \ref{star}.)

 \end{Lemma}

 \begin{proof} 
 Write $A$ for~$\ext(s)$.  
  
 First note that {if} $s=^* s'$, then also $\ext(s) =^* \ext(s')$, so
by \ref{gea.fact} we may  (replacing $s$ by an appropriate  $s' =^*
 s$)
 without loss of generality assume 
 that not only $F(f,s)\preccurlyeq_s^* F(g,s)$ but even  
  $F(f,s)\preccurlyeq_s F(g,s)$. 
  
\medskip
 {\em Part 1.}    
  We first show $\forall k\in A : f(k) \le (h_A\circ g)(k)$.

 So fix $k\in A $.  Let $\eta\in F(f,s)$, 
                             $\eta\trianglelefteq k$.

 Let $\eta'\in  F(f,s)$, where $  \eta \lessdot  \eta'$.

 So $f(k) \le f(\max[\eta]) \le  \min [\eta']$, because $F(s,f)$ gauges~$f$.  
  
 But $g(k)\ge g(\min[\eta]) > \max[\eta]$, as $\eta$ is $g$-weak.  
   By the definition of  
 $\eta \lessdot\eta'$ this means that $h_A(g(k))\ge \min [\eta']$, as  
 $h_A(\cdots) \in A$.  
  
 So $f(k) \le (h_A\circ g)(k)$.   
  
\medskip

 {\em Part 2.} 
 We now consider a general $n\in \N$.  Let $k:= h_A(n)$.  
  Clearly $n\le k\in A=\ext(s)$.  
 So by part 1, $f(n)\le f(k)  \le (h_A\circ g)(k) =  
  (h_A\circ g\circ h_A)(n) $.  
  
So in any case we have $f\le h_A \circ g \circ h_A$, which means
$f\le_A g$. 
 \end{proof}

 \begin{Lemma} \herea123a.\label{not} 
   If $f \prec^*_s g$, then $f \le_{\ext(s)} g$ and $g \not\le_{\ext(s)} f$. 
 \end{Lemma} 
  
 \begin{proof} We already have $f \le_{\ext(s)} g$, so we only 
have to refute $g \le_{\ext(s)} f$. 

 Again replacing $s$ by an appropriate $s' =^* s$ we may assume $f
\prec_s g$. For notational simplicity only we will assume 
that  the *fronts $F(s,f)$ and $F(s,g)$ are actually fronts. 

Let  $A:= \ext(s)$.   Assume $g \le \xx \max(h_A, f)^j $. 
 All except finitely many $\eta\in \int(s)$ have more than $j$ direct 
 successors; find $n$ and $\eta\in F(s,g)\cap S_n$ such that $\eta$ has 
 more than $j$ successors.
 Consider the set $C:=\{ 
 \nu\in F(s,f): \eta \vartriangleleft \nu\}$.  We know that this set is 
 nonempty and even that it it has more than $j$ elements (since each 
 branch through $\eta$ must meet~$C$).    
 We can write $C$ as $C = \{\nu_1, \ldots , \nu_\ell\}$, where  
 $\ell > j$ and $\nu_1 \lessdot \nu_2 \lessdot \cdots \lessdot \nu_\ell$.  
  
  Then $h_A(\min[\nu_i] ) \le \min [ \nu_{i+1}]$ for all~$i$, and 
 also  $f(\min[\nu_i] ) \le \min [ \nu_{i+1}]$, as 
$\nu_i \in F(s,f)$.   Hence  
 $$ \xx \max(h_A, f) ^j (\min [\nu_0])  \le \max[\nu_\ell].$$ 
 But $g(\min[\nu_0] )  > \max[\nu_\ell]$, because $\eta $ is $g$-weak
(as
$\eta\in F(s,g)$).  
 \end{proof}

The last two facts allow us to replace the relation $\le_A$ between
functions {\sl (this is the relation that we are really interested in)\/} by
the relation $\preccurlyeq^*_s$ between the associated fronts 
{\sl (this is
the relation that can be more easily manipulated, by modifying~$s$)}, 
 assuming that $A = \ext(s)$ and that $s$ ``knows enough'' 
about $f$ and~$g$.

 \section{Direct limit}

 We will fix a nonempty
 partially ordered  $(I, {\le})$ in which every countable set
has an upper bound.  
  Later we will consider only the special case $I = \omega_1$.

 \begin{Definition} \herea241. 
 Let $\vec s = (s_i: i\in I)$ be a sequence of zoos with  
 $\forall i < j:  s_j \leqq^* s_i $.    We say that $F$ is a  *front 
 in~$\vec s$ {if} there  is $i\in I$ such that $F \subseteq s_i$ is a
 *front in~$s_i$.   
   (Note:  This implies that each $F\cap s_j$ is a *front in~$s_j$,
 for all $j\ge i$.)   
  
 \end{Definition} 
  
 \begin{Definition}\herea241a. 
    Let $\vec s$ be as above, and let 
 $F$ and $G$ be *fronts in~$s_{i_1}$,~$s_{i_2}$, respectively.  
  
 We will write $F \approx^*_{\vec s} G$ 
  iff there is some $i \ge  i_1, i_2$   such that  we have  
  $F \cap s_i \approx^*_{s_i} G\cap s_i$.  
 Equivalently, we could demand:   
 \begin{quote} For some  $i_* \ge  i_1, i_2$, for all $i\ge i_*$:  
  $F \cap s_i \approx^*_{s_i} G\cap s_i$ 
 \end{quote} 
Clearly, this is an equivalence relation.
  
 Similarly we define $F \preccurlyeq_{\vec s} G$, $F
 \preccurlyeq_{\vec s}^* G$, $F  + n \approx_{\vec s}^* G$,
 etc. (See definition~\ref{prec} and
 fact~\ref{persist}.)  

For example, $F + n \aps G$ iff one or both of the following two
equivalent  conditions hold: 
\begin{enumerate}
\item $\exists i:\  (F\cap s_i ) + n \approx^*_{s_i} (G\cap s_i)$
\item $\exists i_0\forall i\ge i_0:\  (F\cap s_i ) + n \approx^*_{s_i} (G\cap s_i)$
\end{enumerate}

 Given $\vec s$, we write $\F\vec s$ for the set of $\aps$-classes 
 of *fronts.   $\F\vec s$ is naturally partially ordered by~$\prs$. 
  \end{Definition}
  
 \begin{Definition} \herea242. 
 Let $\vec s$ be as above. 
   
 For each~$i$, {if} $s_i = (\S_0^i, \S_1^i, \ldots)$, we let 
  $\wurzel(s_i) $ be 
 the front $\{\wurzel(S_0^i), \wurzel(S_1^i), \ldots \}$.  
 \end{Definition}

 \begin{Definition} \label{suf1}
 $\vec s = (s_i:i\in I)$ is ``sufficiently generic''  if:  
 \begin{enumerate} 
 \item $i < j$ implies $s_j \leqq^* s_i$  
\item For any  $i\in I$, and  any coloring $c$ of $s_i$ (of one of the three
types described in \ref{color.zoo}) there is $j\in I$, $j > i$, such that
$ s_j$ is  almost   $c$-homogeneous. 
 \item For all $f:\N\to \N$ there is $i$ such that $s_i$ is almost 
$f$-strong. 
\item For all $i$ there is $j\ge i$ such that $(\wurzel(s_i)\cap s_j) +1
\preccurlyeq^* \wurzel(s_j)$. 
\end{enumerate} 
\end{Definition}

{\sl 

At the end of this section we will show that (assuming CH) there 
exists a sufficiently generic sequence  $(s_i: i \in \omega_1)$. 

But first we will show  how a sufficiently generic sequence helps
to get the desired clone.    We first show that the set $\F\vec s$ of 
$\approx^ *_{\vec s}$-equivalence classes of *fronts  is linearly ordered
by $\le_{\vec s}$, without last element. 

As a byproduct, we get some more information about this linear order (such as:
every element has a direct successor).  Assuming CH, this information will
be sufficient to characterize this order up to order isomorphism.

We then consider the filter $U$ generated by the sets~$\ext(s)$. The 
results from the previous section will easily show that the map 
$f \mapsto F(s_i,f)$  (for an appropriate  $i = i(f)\in I$) induces 
an isomorphism between the order $\G/U$ and $\F\vec s$.  This 
is enough to prove our  main theorem.

}

\begin{Fact} \label{suf2}
Assume that 
  $\vec s = (s_i:i\in I)$ is sufficiently generic. Then 
\begin{enumerate}\setcounter{enumi}{4}
\item\label{i.many}
 For all~$i$: if $c_0, c_1,\ldots$ are colorings of~$s_i$, then
there is $j>i$ such that $s_j$ is almost homogeneous for each~$c_k$.

\item \label{i.n}
For all *fronts $F,G$ we have $F \prs G$ or $G \prs F$. 
\item \label{i.tri}
For all *fronts $F,G$ with $F
\prs G$, exactly one of the following holds:
\begin{itemize}
\item  $F + n \aps G $ for some (unique) $n$
\item  $F + \infty  \prs G $. 
\end{itemize}
\item \label{i.1}
For all *fronts $F$ there is a *front $G$ such that $F+1 \aps G$
\item \label{i.bound}
Every countable set of *fronts has a $\prs$-upper bound.  Moreover, if
${\mathscr F}_1 $ and ${\mathscr F}_2$ are countable sets of *fronts and
$\forall F_1\in{\mathscr F}_1\, \forall F_2\in{\mathscr F}_2: F_1 \prs F_2$,
then
there is a *front $G$ with
$$\forall F_1\in{\mathscr F}_1\, \forall F_2\in{\mathscr F}_2:
          F_1 \prs G \prs  F_2.$$
\item \label{i.infty}
For all *fronts $F$ there is a *front $G$ such that $F+\infty \prs G$ 
\item \label{i.c}
Whenever  $F+\infty \prs G$, then there are $2^{\aleph_0}$ equivalence classes of fronts between $F $ and~$G$. 
{\small   (In other words: 
 {if} we divide the set $\F(\vec s)$ by the equivalence relation 
 generated by  ``the interval $[x,y]$ is finite'' then
 we get a linear order which is  $2^{\aleph_0}$-dense.) } 
\item \label{i.ext}
For all~$i,j$:  $\ext(s_i)\aps \ext(s_j) $; the equivalence class of these *fronts is the $\prs$-smallest class. 
\item\label{i.pre}
 For all *fronts $G $ not in the class described in
 (\ref{i.ext}), there is a *front $F$ with 
$F + 1 \aps G$. 
 \item \label{i.gauge}
For all $f:\N\to \N$ there is $i$ such that $s_i$ gauges~$f$.  
\end{enumerate}
\end{Fact}

\begin{proof}\

(\ref{i.many})  
We can find a sequence $ i \le i_0 \le i_1 \le \cdots$ in $I$ such
that $s_{i_n} $ is almost $c_n$-homogeneous.  Now let $j$ be any upper
bound of the set $\{i_0,i_1,\ldots\}$, then $s_j$ must be almost
$c_n$-homogeneous for all~$n$.

(\ref{i.n}) By \ref{dichotomy}(1), and \ref{suf1}(1,2).

(\ref{i.tri}) By (\ref{i.many}) and \ref{dichotomy}(3). 

(\ref{i.1}) Easy, using \ref{suf1}(4). 

(\ref{i.bound}) Not needed for our main conclusion, and left to the reader.

(\ref{i.infty}) Use (\ref{i.1}) and (\ref{i.bound}).

(\ref{i.c}) Not needed for our main conclusion, and left to the reader.

(\ref{i.ext}) Clear.

(\ref{i.pre}) Clear.

(\ref{i.gauge})  Use \ref{successor1} and \ref{suf1}(1,2,3).

\end{proof}

The following transfinite construction is now  routine.

 \begin{Conclusion} \herea245. 
 \label{CH1} 
   Assume CH.  Then there is a sufficiently generic sequence 
 $(s_i: i \in\omega_1)$.  
 \end{Conclusion} 
 \begin{proof}

Recall that $\omega_1$ is an uncountable well-ordered set with the
property that for all $j\in \omega_1$ the set $\{i\in \omega_1: i <
j\}$ is countable.  Also recall that each element $i \in \omega_1$ has
a direct successor
$$ i+1 := \min \{j\in \omega_1: i < j \}$$
 Limit points of $\omega_1$
are those elements which are not of the form~$i+1$. The least 
element of $\omega_1$ is called~$0$.

We will use a straighforward bookkeeping argument to take care of 
\ref{suf1}(2--4).

By CH, let $(f_i: i \in \omega_1) $ enumerate all growth functions.
Let $H:\omega_1\times \omega_1\to \omega_1$ be a bijection satisfying
$H(\alpha, \beta)\ge \alpha$ for all $\alpha, \beta\in \omega_1$, and
let $(H_1,H_2)$ be the inverse functions, i.e., $H(H_1(\gamma),
H_2(\gamma) ) = \gamma $ for all $\gamma \in \omega_1$.

Now define a $\leqq^*$-decreasing sequence $(s_i: i \in \omega_1)$ as follows: 

\begin{itemize}
\item $s_0$ is arbitrary. 
\item Assume that all $(s_j: j\le i)$ is already defined.  We will
define~$s_{i+1}$.
Let $\{c^i_j: j \in \omega _1\}$ be the set 
of all coloring functions of~$s_i$. (Again, such an 
enumeration exists because we are assuming CH.)
\\
Now let $s_{i+1}\leqq s_i$ be such that 
\begin{enumerate}
\itm a  $s_{i+1}$ is $f_i$-strong.
\itm b  $(\wurzel(s_i) \cap s_{i+1})+1  \preccurlyeq_{s_{i+1}}
 \wurzel(s_{i+1})$. 
\itm c $s_{i+1}$ is $c^{H_1(i)}_{H_2(i)}$-homogeneous.  (Recall that 
$H_1(i)\le i$, so  $s_{H_1(i)}$ is already defined; 
$c^{H_1(i)}_{H_2(i)}$ is a coloring of~$s_{H_1(i)}$.)
\end{enumerate}
$s_{i+1}$ can easily be obtained in 3 steps $s_i \geqq
s_i^{\rm (a) }\geqq 
s_i^{\rm (b) }\geqq 
s_i^{\rm (c) } = s_{i+1}$, where in each step we satisfy one of the demands
(a), (b), (c). 
  For example, (b) can be realized by using the ``gluing'' step
  from~\ref{examples}.
\item If $j$ is a limit point, then let $s_j$ be any zoo satisfying
$\forall i< j: s_j \leqq^* s_i $. This is possible by \ref{limit2}.  
(See also remark~\ref{finalremark}.)
\end{itemize}

 \end{proof}

 The following facts are easy and well-known:  
  
 \begin{Fact}\herea247. \label{cantor1}\ 
\begin{enumerate}
\item
 Assume CH.  Then there is a unique linear order $(D_1,\le)$  
 with the following properties:  
 \begin{itemize} 
 \item $D_1$ has a smallest element  but no largest element.  
 \item $D_1$ is of cardinality~$\aleph_1$.  
 \item  $D_1$ is $\aleph_1$-dense, i.e.: between any two elements 
 there are uncountably many elements.  
 \item Every countable subset  of~$D_1$ is bounded, and moreover: \\  
For any two countable sets $C,C'\subseteq D_1$  with $C \le C'$
 [i.e.,   $\forall x \in C\, \forall y \in C':  x \le y$] 
there is $c$ with $C \le \{c\} \le  C'$. 
 \end{itemize} 
\item
 Assuming CH, there is also a unique linear order $D_2$ with the following  
 properties:  
 \begin{itemize} 
 \item $D_2$ has a smallest element  but no largest element.  
 \item Every element of~$D_2$ has a direct successor.  
 \item Every element of~$D_2$ (except for the minimal element) has a
 direct predecessor.   
 \item Factoring $D_2$ by the relation 
$$x \sim y  \qquad \Leftrightarrow \qquad \mbox{
 the interval $[x,y]$ is finite}$$
 yields~$D_1$.  
 \end{itemize} 
\end{enumerate}

 \end{Fact}

 \begin{proof}  
  A back-and-forth argument, similar to Cantor's theorem characterizing 
  the rationals as the unique dense linear order without endpoints. 
  
 $D_2$ can be obtained as the lexicographic order on  
 $(\{\min D_1\} \times \N) \cup (D_1\setminus\{\min D_1\}) \times \Z $.

 \end{proof} 
  
 \begin{Conclusion}\herea248. \label{done}
 Assume CH.  
  
 Let $\vec s = (s_i: i\in \omega_1)$ be sufficiently generic, and let 
 $U$ be the  filter generated by $(\ext( s_i): i \in \omega_1)$. 
 Then 
 \begin{itemize} 
 \itm a $U$ is an ultrafilter. 
 \itm b $(\F\vec s , {\prs}) $ is order isomorphic to~$D_2$. 
 \itm c 
  The set $\G/ {\sim_U}$, ordered by~$\preccurlyeq_U$, is order
  isomorphic to~$D_2$.   
 \itm d 
  Letting $\langle U \rangle$ be the clone generated 
by $\C_\id \cup \{h_A:A\in U\}$, the 
 interval $[\langle U\rangle, \O] $
 in the clone lattice is order isomorphic 
 to the Dedekind completion of~$D_2$.  
 \end{itemize} 
 \end{Conclusion} 
 \begin{proof}  (a) is not needed and left to the reader.  

\bigskip

(b) By \ref{suf2} and \ref{cantor1}.

 \bigskip
(c)
We define a map $K$ from $\G$ to the set of *fronts in $\vec s$:\\ For any
$f\in \G$, pick $i\in I$ such that $s_i$ is almost $f$-strong and gauges~$f$. 
Let $K(f) = F(s_i,f)$.   (Recall  that for all
 $j\ge i$, $F(s_i,f) \approx^*_{s_j} K(f) \cap s_j$.)

We now claim that 
\begin{itemize}
\item[$(*)$] $ f \le_U g$ iff $K(f) \prs K(g)$.
\end{itemize}

To prove this claim, fix $f$ and~$g$.  Pick some sufficiently large
  $j\in\omega_1$.  
By \ref{suf2}, one of the following cases holds:

 \begin{itemize}
\item[(i)] $K(f) \aps K(g)$, so  $F(s_j,f) \approx_{s_{j}}^ * F(s_j,g)$
\item[(ii)] $K(f) + 1 \prs K(g)$, so $F(s_j,f) \prec_{s_{j}}^ * F(s_j,g)$
\item[(iii)] $K(g) + 1 \prs K(f)$, so $F(s_j,g) \prec_{s_{j}}^ * F(s_j,f)$
 \end{itemize}


In the first case, 
 \ref{leq} implies $ f\le_{\ext(s_{j})} g$, hence $f\le_U g$. 

In the second case, we have $F(s_i,f) \prec_{s_i}^* F(s_i,g) $ for all $i\ge j$, 
so by~\ref{not}
$$ \forall i\ge j: \qquad    f \le_{\ext(s_i)} g \quad \mbox{ and } 
\quad g \not\le_{\ext(s_i)} f$$
which implies $f\le_U g$, and $g\not\le_U f$. 

In the third case we get similarly $g\le_U f$, and $f\not\le_U g$. 

So in each case the desired equivalence holds, and $(*)$ is proved.

 Hence $K$ induces an order isomorphism between $\G/U$ and $\F\vec s $. 

\bigskip

(d) follows from our discussion in section~\ref{section.reduce}.  
 \end{proof}

\begin{Remark}\label{finalremark}
It is clear that the full strength of CH is not necessary for this
construction.    Martin's axiom MA (even the version for $\sigma$-centered
forcing notions) is easily shown to imply  an 
analogue of corollary~\ref{limit2}, 
in which ``countable'' is replaced by ``of size $< 2^{\aleph_0}$''.  
%
This allows us to modify the construction in \ref{CH1} to 
a transfinite induction of length $2^{\aleph_0}$, which shows 
that already MA implies the existence of a sufficiently 
generic sequence. 

Thus, the conclusion of our theorem is also consistent with 
the negation of CH.

\end{Remark}

 \bibliographystyle{plain} \bibliography{other,listb,listx}

 \end{document}